\documentclass[a4paper,english]{article}

\usepackage[utf8]{inputenc}
\usepackage{mathrsfs}
\usepackage{babel}
\usepackage{amsmath,amsthm,amssymb,amsxtra}
\usepackage{colortbl}
\usepackage{hyperref}
\usepackage{enumitem}
\usepackage{pdfsync}
\usepackage{colortbl}

\newcommand*{\scrF}{\ensuremath{\mathscr{F}}} 
\newcommand*{\scrN}{\ensuremath{\mathscr{N}}} 
\newcommand{\cqd}{\hfill$\Box$}
\newcommand*{\caB}{\ensuremath{\mathcal{B}}}	
\newcommand*{\caC}{\ensuremath{\mathcal{C}}}	
\newcommand*{\caF}{\ensuremath{\mathcal{F}}}	
\newcommand*{\caH}{\ensuremath{\mathcal{H}}}	
\newcommand*{\caN}{\ensuremath{\mathcal{N}}}	
\newcommand*{\caS}{\ensuremath{\mathcal{S}}}	
\newcommand*{\D}{\mathbb{D}}									

\newcommand*{\N}{\mathbb{N}}									
\newcommand*{\R}{\mathbb{R}}									
\newcommand*{\Rd}{{\mathbb{R}^d}}							

\newcommand*{\eps}{\varepsilon}								

\newcommand*{\E}{\mathbb{E}}									
\renewcommand*{\P}{\mathbb{P}}								

\newcommand*{\e}{\mathrm{e}}
\newcommand*{\ii}{\mathrm{i}}
\newcommand*{\intd}{\mathrm{d}}								
\newcommand{\Infkt}[1]{1_{ #1 }}							

\DeclareMathOperator*{\supp}{supp}						
\DeclareMathOperator*{\dist}{dist}						

\numberwithin{equation}{section}
\newtheorem{lem}{Lemma}[section]
\newtheorem{thm}[lem]{Theorem}
\newtheorem{prop}[lem]{Proposition}

\newtheorem{numrem}[lem]{Remark}

\allowdisplaybreaks

\begin{document}

\title{Logarithmic asymptotics of the densities of SPDEs driven by spatially correlated noise}
\author{Marta Sanz-Sol{\'e}\thanks{Supported by the grant MTM 2012-31192 from the Direcci\'on General de Investigaci\'on, Ministerio de Economia y Competitividad, Spain.} \ and Andr{\'e} S{\"u}\ss\thanks{Supported by the grant MTM 2012-31192 from the Direcci\'on General de Investigaci\'on, Ministerio de Economia y Competitividad, Spain.}}
\date{\today}

\maketitle

\begin{abstract}
We consider the family of  stochastic partial differential equations indexed by a parameter $\eps\in(0,1]$,
\begin{equation*}
Lu^{\eps}(t,x) = \eps\sigma(u^\eps(t,x))\dot{F}(t,x)+b(u^\eps(t,x)),
\end{equation*}
$(t,x)\in(0,T]\times\Rd$ with suitable initial conditions. In this equation, $L$ is a second-order partial differential operator with constant coefficients, $\sigma$ and $b$ are smooth functions and $\dot{F}$ is a Gaussian noise, white in time and with a stationary correlation in space. Let $p^\eps_{t,x}$ denote the density of the law of $u^\eps(t,x)$ at a fixed point $(t,x)\in(0,T]\times\Rd$. We study the existence of $\lim_{\eps\downarrow 0} \eps^2\log p^\eps_{t,x}(y)$ for a fixed $y\in\R$. The results apply to classes of stochastic wave equations with $d\in\{1,2,3\}$ and stochastic heat equations with $d\ge1$.
\medskip

\medskip

\end{abstract}

\section{Introduction}
In this article, we consider the family of stochastic partial differential equations (SPDEs) indexed by a parameter $\eps\in(0,1]$ defined by
\begin{equation}
	Lu^{\eps}(t,x) = \eps\sigma(u^\eps(t,x))\dot{F}(t,x) + b(u^\eps(t,x)),
	\label{eq:SPDE}
\end{equation}
$(t,x)\in(0,T]\times\Rd$, $d\ge 1$, with suitable initial conditions. Here $L$ is a second-order partial differential operator, typical examples are the wave and the heat operators; $\sigma,b:\R\rightarrow\Rd$ are smooth functions; $\dot{F}$ is a Gaussian noise, white in time and with a stationary correlation in space.

Equation \eqref{eq:SPDE} describes a nonhomogeneous initial value problem subject to nonlinear small random fluctuations. The results of this paper are a contribution to the study of the behavior of \eqref{eq:SPDE} as $\eps\downarrow 0$ and therefore, when the random perturbations disappear. More precisely, denote by $p^\eps_{t,x}$ the density of the random variable $u^{\eps}(t,x)$ at a given point $(t,x)\in(0,T]\times\Rd$. We will determine the set of points $y\in\R$ for which one can derive upper and lower bounds for $\lim_{\eps\downarrow 0} \eps^2\log p^\eps_{t,x}(y)$. We will identify these bounds and refer to them as {\it Varadhan estimates} or {\it logarithmic estimates}. We will consider examples of stochastic wave equations with $d=\{1,2,3\}$ and stochastic heat equations with $d\ge 1$.

For solutions to stochastic differential equations driven by a standard Brownian motion, $\{X_t, t\ge 0\}$, this question is equivalent to the analysis of the density of $X_t$, when $t\downarrow 0$. Under ellipticity conditions and with analytical methods, it has been firstly studied in \cite{varadhan-1-67, varadhan-2-67}. Using Malliavin calculus and large deviation estimates, Varadhan's results have been extended in \cite{leandre-2-87,leandre-1-87} under hypoelliptic assumptions.

The method of \cite{leandre-2-87,leandre-1-87} has been applied in \cite{leandrerusso} to establish Varadhan estimates for an example of hyperbolic SPDE: an It\^o equation with two-dimensional parameter. In \cite[Propositions 4.4.1, 4.4.2]{nualartstflour}, a general formulation of that method is given, providing a systematic approach to the study of Varadhan estimates for families of Wiener functionals subject to small perturbations of their sample paths. For example, it has been used in  \cite{K-M-SS-02} to extend the results of \cite{leandrerusso}, and in \cite{milletsanz1996} for a stochastic heat equation with boundary conditions.


Similarly as in \cite{marquezsarra}, the aim of this paper is to study Varadhan estimates for the class of SPDEs defined by \eqref{eq:SPDE}. However, in comparison with this reference, there are two additional substantial contributions in our results. Firstly, the scope of application of the theory presented in this article is larger. Indeed, we are able to deal with cases where the fundamental solution corresponding to the operator $L$ is a measure, like for example, the stochastic wave equation in spatial dimension $d=3$. Secondly, in \cite[Theorem 1.2]{marquezsarra} it is not clear for what values of $y\in\R$ the claim $\lim_{\eps\downarrow 0} \eps^2\log p^\eps_{t,x}(y)=-I(y)$, where $I$ is the rate function, holds. This statement requires $p^\eps_{t,x}(y)>0$ for $\eps$ small enough, but this problem is not discussed in \cite{marquezsarra}. Also in \cite[Proposition 5.1]{marquezsarra}, it is assumed that the interior of the topological support of the law of $u^\eps(t,x)$ is described in a way that we do 
not see justified. In this paper, these issues are rigorously  addressed.

We now describe the contents of this article. In Section \ref{sec:2}, we formulate the basic assumptions used throughout the paper, we give a rigorous formulation of \eqref{eq:SPDE} and quote two fundamental results concerning the existence of a unique {\it random field solution} to \eqref{eq:SPDE}, and on the {\it existence and smoothness} of the density $p^\eps_{t,x}$ (see Theorems \ref{basic}, \ref{density}, respectively). In Theorem \ref{thm:LDPdensity} we state the main result of the paper on the logarithmic estimates.

Section \ref{sec:3} is devoted to the proof of Theorem \ref{thm:LDPdensity}. To obtain the upper bound, we check that
$u^\eps(t,x)$ is Malliavin differentiable of any order, and that the corresponding Malliavin-Watanabe norm is uniformly bounded in $(t,x)$ and $\eps$. We also prove a quantitative result on the dependence on $\eps$ of the $L^p$ norm of the inverse of the Malliavin matrix corresponding to $u^\eps(t,x)$. Notice that in Theorem \ref{thm:LDPdensity}, the upper bound still makes sense if $\lim_{\eps\downarrow 0} \eps^2\log p^\eps_{t,x}(y)=-\infty$.

To establish the lower bound, we prove that the mapping $\eps\mapsto u(t,x; \omega+\eps^{-1}h)$, where $h$ is an admissible shift for the space of paths $\Omega$, is differentiable in the $\mathbb{D}^\infty$ topology of Malliavin calculus, and that the mapping given in \eqref{eq:defbhphi} is onto. Then,
in order to give full meaning to the lower bound \eqref{eq:ldpdenslower}, it is relevant to know for which set of $y\in\R$, $p^\eps_{t,x}(y)$ is strictly positive for $\eps$ small enough, and whether the function $I$ is finite. In the analysis of these questions, the characterization of the topological support of the law of the random variables $u^\eps(t,x)$, $\eps\in(0,1]$ plays a crucial role. Each one of these random variables are a nonlinear functional $\Phi$ (not depending on $\eps)$ of the driving Gaussian noise $\eps F$. Hence, one should expect the support to be independent of $\eps$. We postpone the proof of a characterization of the support of $u^\eps(t,x)$, which in particular shows its independence of $\eps$, to Section \ref{sec:A3}.

The regularity (in the Malliavin sense) of $u^\eps(t,x)$ established in Lemma \ref{lem:nondegenerateI}, combined with \cite[Proposition 4.1.1, 4.1.2]{nualartstflour} imply that the support of $u^\eps(t,x)$ is a nonempty closed interval and that $p^\eps_{t,x}(y)>0$ for all $y$ in the interior of that set. We also prove in Proposition \ref{finite} that, in these points, $I(y)<\infty$, and also that if the function $b$ is bounded then  $\{y\in\R: I(y)<\infty\}= \R$ (see Proposition \ref{prop:positivitydensity}).

Section \ref{sec:A3} is devoted to the characterization of the topological support of the law of $u^\eps(t,x)$ (see Theorem \ref{ta.1}). The relevant reference is \cite{delgado--sanz-sole012}, where a characterization of the support of the law of a stochastic wave equation in spatial dimension $d=3$ with vanishing initial conditions in H\"older norm is established. In comparison with that work, here the SPDE is more general but, instead of considering the sample paths of the solution to \eqref{2.2}, we take its value at a fixed point $(t,x)$. This makes the analysis significantly easier.

In Section \ref{sec:applications}, we give two examples where the main result is applied: a class of stochastic wave equations with $d\in\{1,2,3\}$ and a class of stochastic heat equation with $d\ge 1$. For the former, owing to results on large deviations, we have $I=J$ and therefore the equality between the upper and lower bounds.

Throughout the paper, we have to deal with different types of evolution equations, including some classes of Hilbert space-valued equations. To provide the suitable background, we prove in the Appendix a result on the existence and uniqueness of random field solution for a very general class of SPDEs.


\section{Preliminaries and statement of the main result}\label{sec:2}

Let $\caC_0^\infty(\R_+\times\Rd))$ denote the space of infinitely differentiable functions defined on $\R_+\times\Rd$ with compact support. On a given probability space $(\Omega,\scrF,\P)$, we consider a Gaussian stochastic process
$F=(F(\phi);\; \phi\in\caC_0^\infty(\R_+\times\Rd))$ with mean zero and covariance functional
\[ J(\phi,\psi) := \E[F(\phi)F(\psi)] = \int_0^\infty\int_\Rd \big(\phi(t)\star\tilde{\psi}(t)\big)(x)\Gamma(\intd x)\intd t, \]
where $\tilde{\psi}(t,x) := \psi(t,-x)$, the symbol ``$\star$" denotes the convolution operator on $\R^d$, and $\Gamma$ is a nonnegative, nonnegative definite, tempered measure on $\Rd$. We know by \cite[Chapter VII, Th\'{e}or\`{e}me XVIII]{schwartz} that there exists a nonnegative tempered measure $\mu$ on $\Rd$ such that $\caF\mu = \Gamma$, where $\caF$ denotes the Fourier transform operator given by
\[ \caF\phi(\xi) = \int_{\R^d} \phi(x) \e^{-2\pi\ii\langle\xi,x\rangle}\intd x. \]

Following \cite{dalangfrangos}, the process $F$ can be extended to a worthy martingale measure
$M=(M_t(A);\; t\in\R_+, A\in\caB_b(\Rd))$ where $\caB_b(\Rd)$ denotes the bounded Borel subsets of $\Rd$. This is achieved by approximating indicator functions $\Infkt{A}$, $A\in\caB_b(\R_+\times\Rd)$ by functions in $\caC^\infty_0(\R_+\times\Rd)$, and thus extending the functional $\phi\mapsto F(\phi)$ to an $L^2(\Omega)$-valued measure $A\mapsto F(\Infkt{A})$. Then we define
\[ M_t(A) := F(\Infkt{[0,t]\times A}), \]
for all $t\in\R_+$ and $A\in\caB_b(\Rd)$.

 Throughout this article we use the filtration
\[ \scrF_t := \sigma\left(M_s(A);\; s\in[0,t], A\in\caB_b(\Rd)\right)\vee\scrN, \]
$t\in\R_+$, where $\scrN$ is the $\sigma$-field generated by the $\P$-null sets.

The SPDE \eqref{eq:SPDE} is expressed in the {\it mild formulation}, as follows,
\begin{align}
	u^\eps(t,x) = & w(t,x) + \eps\int_0^t\int_\Rd \Lambda(t-s,x-z)\sigma(u^\eps(s,z))M(\intd s,\intd z) \notag\\
								& + \int_0^t\int_\Rd \Lambda(t-s,x-z)b(u^\eps(s,z))\intd z\intd s,	\label{2.2}
\end{align}
$(t,x)\in(0,T]\times \Rd$, where $\Lambda$ denotes the fundamental solution to the associated PDE, $Lu=0$, and $w$ is the contribution of the initial conditions. For $\eps=1$, we will write $u(t,x)$ instead of $u^1(t,x)$.

We will consider the following assumptions:

\begin{enumerate}[label=\bfseries (A\arabic{enumi}),ref=\bfseries (A\arabic{enumi})]
  \item\label{itm:finiteLambda} The mapping $t\mapsto\Lambda(t)$ is a deterministic function with values in the space of non-negative tempered distributions with rapid decrease such that
  \[ \int_0^T\int_\Rd \big(\Lambda(s)\star\tilde{\Lambda}(s)\big)(x)\Gamma(\intd x)\intd s = \int_0^T \int_\Rd |\caF\Lambda(s)(\xi)|^2\mu(\intd\xi)\intd s < \infty. \]
  Moreover, for all $t\in(0,T]$, $\Lambda(t)$ is a nonnegative measure, and there exists $\delta>0$ such that
  \begin{equation}\label{eq:conditiondelta}
  	\int_0^t\Lambda(s)(\Rd) ds \le C t^\delta.
  \end{equation}
  \item\label{itm:conditionic} The mapping $(t,x)\mapsto w(t,x)$ is deterministic, continuous and
  \[ \sup_{{(t,x)\in[0,T]\times\Rd}} |w(t,x)| < \infty. \]
\end{enumerate}

\begin{numrem}
\label{general}
Later on, we will refer to  \cite{dalang} and also to \cite{dalangquer} for results on the stochastic integral in \eqref{2.2}, and on the existence and uniqueness of solution. These are proved assuming that $\sup_{t\in[0,T]} \Lambda(t)(\Rd)<\infty$.
It can be easily checked that they also hold assuming \eqref{eq:conditiondelta}.
\end{numrem}

Throughout the paper the following notation will be used. Let $\Lambda$ be as in hypothesis \ref{itm:finiteLambda}. For any $s\in[0,T]$, set
\begin{align*}
	J_1(s) & := \int_\Rd\big(\Lambda(s)\star\tilde{\Lambda}(s)\big)(z)\Gamma(\intd z) = \int_\Rd |\caF\Lambda(s)(\xi)|^2\mu(\intd\xi), \\
	J_2(s) & := \Lambda(s)(\Rd),\\
	g_1(t) & := \int_0^t J_1(s)\intd s. 
\end{align*}
Notice that \ref{itm:finiteLambda} implies $g_1(T) < \infty$.

In \eqref{2.2}, the last integral denotes the convolution $\int_0^t (\Lambda(t-s)\star b(u(s,\cdot))(x)\intd s$ in the space variable, defined pathwise.
As for the stochastic integral (also termed {\it stochastic convolution}), we refer to the construction given in \cite{dalang} (see \cite[Section 2.3]{dalangquer} for a summary).

Let $\mathcal{S}(\Rd)$ be the set of Schwartz functions and denote by $\mathcal{H}$ the Hilbert space obtained by completion of the set $\mathcal{S}(\R^d)$ with the inner product
\[ \langle \phi,\psi\rangle=\int_{\R^d}(\phi\star\psi)(x) \Gamma(\intd x) = \int_{\R^d}\caF\phi(\xi)\overline{\caF\psi(\xi)} \mu(\intd\xi). \]
Set ${\caH_T} := L^2([0,T];\caH)$. The Gaussian process $F$ can be extended to an isonormal process $F=(F(\phi);\; \phi\in{\caH_T})$ in the sense of \cite[Definition 1.1.1]{nualart}.

It is useful to identify the isonormal process $F$ with a $\caH$-valued cylindrical Wiener process. As shown in \cite{dalangfrangos}, by an approximation procedure we define  $W_t(\phi)=F(1_{[0,t]}\phi)$, $t\in[0,T]$, $\phi\in\caH$. Consider a complete orthonormal system (CONS) of $\caH$ denoted by $(e_k)_{k\in\N}$. Then,
\[ W=\{W^k(t):=W_t(e_k), t\in[0,T], k\in\N\} \]
defines a sequence of independent standard Brownian motions. Conversely, the process $(F(\phi) = \sum_{k\in\N}\int_0^T \langle\phi(t),e_k\rangle_{\caH}\intd W^k(t), \phi\in\caH_T)$ is an isonormal Gaussian process.

As has been established in \cite{dalangquer}, there is an equivalence between the stochastic integral in the sense of \cite{dalang} and the stochastic integral with respect to the cylindrical Wiener process $W$ (see e.g. \cite{dapratozabczyk}). In particular, the stochastic integral in \eqref{2.2} is equal to
\begin{equation*}
\sum_{k\in\N}\int_0^t \langle \Lambda(t-s,x-\ast) \sigma(u^\eps(s,\ast)),e_k\rangle_{\caH} \intd W^k(s).
\end{equation*}

Appealing to \cite[Theorem 4.3]{dalangquer} and to Remark \ref{general}, for any fixed $\eps\in(0,1]$, there exists a stochastic process
$\{u^\eps(t,x), (t,x)\in[0,T]\times \Rd\}$ such that \eqref{2.2} holds for any $(t,x)\in[0,T]\times \Rd$ a.s. This is termed a {\it random-field solution} to \eqref{2.2}. More precisely, we have the following result.

\begin{thm}
\label{basic}
If Hypotheses \ref{itm:finiteLambda} and \ref{itm:conditionic} are satisfied and $\sigma$ and $b$ are Lipschitz continuous functions, then \eqref{2.2} has a unique random-field solution. Among other properties, this solution is $L^2$-continuous and for any $p\in[1,\infty)$
\[ \sup_{\eps\in(0,1]}\sup_{(t,x)\in[0,T]\times\Rd} \E\big[|u^\eps(t,x)|^p\big] < \infty. \]
\end{thm}


We are interested in the family of densities of the probability law of the solution $u^\eps(t,x)$, $\eps\in(0,1]$ at every fixed point
$(t,x)\in(0,T]\times\Rd$. For this, we describe the abstract Wiener space that will be used as framework for the application
of the Malliavin Calculus (see \cite{nualart}).

Let $(\bar\Omega,\bar{\mathcal{G}}, \bar \mu)$ be the canonical space of a standard real-valued Brownian motion on $[0,T]$. With the equivalence shown before, we can identify the canonical probability space of the stochastic process $F$ with that corresponding to a sequence of independent standard Brownian motions $(\Omega,\mathcal{G},\mathbb{P})=(\bar\Omega^{\N}, \bar{\mathcal{G}}^{\otimes\N}, \bar{\mu}^{\otimes\N})$. This will be the underlying probability space in this article.

Consider the Hilbert space $H$ consisting of sequences $(h^k)_{k\in\N}$ of functions $h^k:[0,T]\rightarrow \R$ which are absolutely continuous with respect to the Lebesgue measure and such that $\|h\|_{H}^2 = \sum_{k\in\N}\int_0^T|\dot h^k(s)|^2 ds < \infty$, where $\dot h^k$ refers to the derivative of $h^k$ defined almost everywhere. There is an isometry $I: H \to\caH_T$ defined by $i(h)(t)= \sum_{k\in\N}\dot{h}^k(t)e_k$. In the sequel we will identify the Hilbert spaces $H$ and $\caH_T$ and by an abuse of notation, we will write $i(h)=h$.
The triple $(\Omega,H,\mathbb{P})$ is the abstract Wiener space that we shall use as framework for the Malliavin calculus.
\medskip

Let us introduce some additional assumptions:
\begin{enumerate}[label=\bfseries (A\arabic{enumi}),ref=\bfseries (A\arabic{enumi})]
\setcounter{enumi}{2}
  \item\label{itm:boundsLambda} There exist positive constants $C,\gamma>0$ and $t_0\in(0,T]$ such that for all $t\in[0,t_0]$,
  \[ Ct^{\gamma} \leq  \int_0^t J_1(s) \intd s = g_1(t). \]
  \item\label{itm:sigmabCinf} The functions $\sigma$ and $b$ are infinitely differentiable with bounded derivatives of any order greater or equal than one.
  \item\label{itm:sigmabound} The function $\sigma$ satisfies $\inf_{x\in\R} |\sigma(x)|=\sigma_0>0$.
\end{enumerate}

The following result in \cite{nualartquer} establishes the existence and regularity of the densities for the solution to \eqref{2.2} at any
point $(t,x)\in(0,T]\times\Rd$.

\begin{thm}
\label{density}
  Fix $(t,x)\in(0,T]\times\Rd$ and $\eps\in(0,1]$. Assume \ref{itm:finiteLambda}, \ref{itm:conditionic}, \ref{itm:boundsLambda}, \ref{itm:sigmabCinf} and \ref{itm:sigmabound}. Then the law of $u^\eps(t,x)$ is absolutely continuous with respect to the Lebesgue measure on $\R$ and its density, denoted by $p^\eps_{t,x}$, is an infinitely differentiable function.
\end{thm}

The last relevant assumption is the following.

\begin{enumerate}[label=\bfseries (A\arabic{enumi}),ref=\bfseries (A\arabic{enumi})]
\setcounter{enumi}{5}
  \item\label{itm:LDPueps} For every $(t,x)\in(0,T]\times\Rd$ the family $(u^\eps(t,x))_{\eps\in(0,1]}$ satisfies a large deviation principle on $\R$ with rate function $J$.
\end{enumerate}
We refer the reader to \cite{dembozeitouni93} for notions and results on large deviations.

 We are now in a position to formulate the main result of this paper. It is about the behaviour of the density $p^\eps_{t,x}(y)$ at every fixed
 $(t,x)\in(0,T]\times\Rd$ and $y\in\R$, as $\eps\rightarrow0$. It will be proved by using the method introduced in \cite{leandre-2-87,leandre-1-87} (see \cite{nualartstflour} for a general formulation).

\begin{thm}\label{thm:LDPdensity}
\begin{enumerate}[label=(\roman{enumi}),ref=(\roman{enumi})]
	\item\label{itm:LPDdensityi} Fix $(t,x)\in(0,T]\times\Rd$ and assume \ref{itm:finiteLambda}, \ref{itm:conditionic}, \ref{itm:boundsLambda}, \ref{itm:sigmabCinf}, \ref{itm:sigmabound} and \ref{itm:LDPueps}. Then for any $y\in \R$,
\begin{equation}
	\lim_{\eps\downarrow0} \eps^2\log p^\eps_{t,x}(y) \leq -J(y).
	\label{eq:ldpdensupper}
\end{equation}
	\item\label{itm:LPDdensityii} Let $(t,x)\in(0,T]\times\Rd$. Assume \ref{itm:finiteLambda}, \ref{itm:conditionic}, \ref{itm:boundsLambda}, \ref{itm:sigmabCinf} and \ref{itm:sigmabound}. Fix $y\in\R$ in the interior of the topological support of the law of $u(t,x)$. Then,
\begin{equation}
	\lim_{\eps\downarrow0} \eps^2\log p^\eps_{t,x}(y) \geq -I(y),
	\label{eq:ldpdenslower}
\end{equation}
with
\begin{equation}
  I(y) = \inf \bigg\{\frac{1}{2}\|h\|_{\caH_T}^2; \; h\in \caH_T, \Phi_{t,x}^{h} = y \bigg\},
  \label{eq:ratefunc}
\end{equation}
and where $\Phi_{t,x}^{h}\in\R$ is defined by
\begin{align}
  \Phi_{t,x}^{h}
  = w(t,x) & + \big\langle\Lambda(t-\cdot,x-\ast)\sigma(\Phi^{h}_{\cdot,\ast}),h\big\rangle_{\caH_T} \notag\\
  & + \int_0^t\int_\Rd \Lambda(t-s,x-z) b(\Phi^{h}_{s,z})\intd z\intd s.
\label{eq:defbhphi}
\end{align}
\end{enumerate}
\end{thm}

We end this section with some important comments on these statements. The existence and uniqueness of a solution to \eqref{eq:defbhphi} follows from Theorem \ref{thm:existenceanduniqueness} in the Appendix.
Theorem  \ref{thm:LDPdensity} makes sense for those $y\in\R$ such that $p^\eps_{t,x}(y)>0$ for all $\eps$ sufficiently small, and the lower bound in \eqref{eq:ldpdenslower} is nontrivial only if $I(y)<+\infty$.
In the last part of Section 3.2, we show the connection between these properties and the topological support of the law of $u^\eps(t,x)$.

Under some additional assumptions, in Section \ref{sec:A3} we will prove a characterization of the topological support of $u^\eps(t,x)$, $\caS$, that exhibits its independence on $\eps$. Proposition \ref{prop:positivitydensity} shows that if $b$ is bounded, $\caS = \mathbb{R}$.

We prove in Proposition \ref{prop:positivitydensity} that $I(y)$ in \eqref{eq:ratefunc} is finite for any $y$ in the interior of $\mathcal{S}$. This uses the characterization of the support. 

Assume that $y$ belongs to the interior of $\mathcal{S}$. Then,  
 \cite[Proposition 4.1.2]{nualartstflour} yields that  $p^\eps_{t,x}(y)>0$ for all $\eps\in(0,1]$.

In Section \ref{sec:applications} we show that for a class of stochastic wave equations, the hypotheses of Theorem \ref{thm:LDPdensity} are satisfied, and $J$ and $I$ in \eqref{eq:ldpdensupper} and \eqref{eq:ldpdenslower} respectively,  are identical. Hence, for any $y$ in the interior of the support of the law we have
\begin{equation*}
	\lim_{\eps\downarrow0} \eps^2\log p^\eps_{t,x}(y) = -I(y),
\end{equation*}
with $I$ defined in \eqref{eq:ratefunc}. With some restrictions, Theorem \ref{thm:LDPdensity} also applies to the stochastic heat equation.

Throughout this article we use the notation $C$ for generic constants that may change from one expression to another. As for the
notations and notions of Malliavin calculus, we refer to \cite{nualartstflour,nualart}.


\section{Proof of the main result}\label{sec:3}
The two parts of Theorem \ref{thm:LDPdensity} will be established separately, applying the methods introduced in \cite{leandre-2-87,leandre-1-87} and extended to an abstract setting in
\cite{nualartstflour} (see Propositions \ref{prop:bound442}, \ref{prop:bound441} below).

By $\Vert \cdot\Vert_{k,p}$, $k\in\N$, $p\in[1,\infty)$, we will denote the norm in the space $\D^{k,p}$ (the Watanabe-Sobolev spaces), and by $\Vert \cdot\Vert_{p}$, the $L^p(\Omega)$ norm. We say that a random variable $X:\Omega\rightarrow\R$ is non-degenerate if $X\in\D^{\infty}=\bigcap_{k\in\N}\bigcap_{p\in[1,\infty)} \D^{k,p}$ and the random variable $\gamma_X := \Vert DX\Vert^2_{\caH_T}$ satisfies $\gamma_X^{-1}\in\bigcap_{p\in[1,\infty)} L^p(\Omega)$, were $D$ denotes the Malliavin derivative. The law of a non-degenerate random variable possesses an infinitely differentiable density.

For the upper bound \eqref{eq:ldpdensupper}, we rely on the following result.

\begin{prop}[\cite{nualartstflour}, Proposition 4.4.2]\label{prop:bound442}
Let $(F^\eps)_{\eps\in(0,1]}$ be a family of non-degenerate random variables satisfying
\begin{enumerate}[label=(\roman{enumi}),ref=(\roman{enumi})]
	\item\label{itm:i} $\sup_{\eps\in(0,1]} \|F^\eps\|_{k,p}<\infty$ for all $k\in\N$ and $p\in[1,\infty)$;
	\item\label{itm:ii} for any $p\in[1,\infty)$ there exists $N_p\in[1,\infty)$ such that $\|\gamma_{F^\eps}^{-1}\|_p\leq \eps^{-N_p}$;
	\item $(F^\eps)_{\eps\in(0,1]}$ obeys a large deviation principle on $\R$ with rate function $J$.
\end{enumerate}
Then
\[ \limsup_{\eps\downarrow0} \eps^2\log p^\eps(y) \leq -J(y), \]
where $p^\eps$ denotes the density of $F^\eps$.
\end{prop}

We denote by $\caC^1(\caH_T;\R)$ the set of all Fr{\'e}chet differentiable real functions $F$ defined on $\caH_T$. For such deterministic functions, we shall use the notation $\bar{D}F$ for its Fr\'echet derivative and set $\bar{\gamma}_F=\Vert \bar{D}F\Vert^2_{\caH_T}$.

 The lower bound \eqref{eq:ldpdenslower} will be established using the following Proposition.

\begin{prop}[\cite{nualartstflour}, Proposition 4.4.1]\label{prop:bound441}
Let $(F^\eps)_{\eps\in(0,1]}$ be a family of non-degenerate random variables. Let $\psi\in\caC^1(\caH_T;\R)$ be such that for all $h\in \caH_T$
\begin{equation}
	\lim_{\eps\downarrow0}\frac{F^\eps(\omega+\eps^{-1}h) - \psi(h)}{\eps} = \caN(h)
	\label{eq:prop441}
\end{equation}
in the $\D^\infty$-topology, where $\caN$ is a random variable belonging to the first Wiener chaos with variance $\bar{\gamma}_\psi(h) = \|\bar{D}\psi(h)\|^2_{\caH_T}$. Then
\begin{align*}
	\liminf_{\eps\downarrow0} \eps^2\log p^\eps(y)
	& \geq -\frac{1}{2}d_R^2(y) \notag\\
	& := -\frac{1}{2}\inf \big\{ \|h\|_{\caH_T}^2; \; h\in \caH_T, \psi(h) = y, \bar{\gamma}_\psi(h)>0 \big\},
\end{align*}
where $p^\eps$ denotes the density of $F^\eps$.
\end{prop}
Let us point out that in the proof of this proposition it is implicitly assumed that $p^\eps(y)>0$ for any $0<\eps<\eps_0$, with $\eps_0$ small enough.


\subsection{Upper bound}

The objective of this section is to apply Proposition \ref{prop:bound442} to the family of random variables  $F^\eps:=u^\eps(t,x)$, $\eps\in(0,1]$, given in \eqref{2.2}, where $(t,x)\in(0,T]\times\Rd$ is fixed.
We will assume that \ref{itm:LDPueps} holds and check that \ref{itm:finiteLambda}--\ref{itm:sigmabound} imply the validity of \ref{itm:i} and \ref{itm:ii} in Proposition \ref{prop:bound442}. This will prove the statement \ref{itm:i} in Theorem \ref{thm:LDPdensity}.

\begin{lem}\label{lem:nondegenerateI}
  Under the conditions \ref{itm:finiteLambda}, \ref{itm:conditionic} and \ref{itm:sigmabCinf} the Assumption \ref{itm:i} in Proposition \ref{prop:bound442} holds. More precisely, we have
  \[ \sup_{{\eps\in(0,1]}}\sup_{{(t,x)\in[0,T]\times\Rd}} \E\Big[\big\|D^ku^\eps(t,x)\big\|_{{\caH_T}^{\otimes k}}^p\Big] < \infty. \]
\begin{proof}
This follows along the same lines as in \cite[Proposition 6.1]{nualartquer}. The difference is that here we are considering a family of SPDEs depending on a parameter $\epsilon\in(0,1]$, and obtain that the norm is uniformly bounded in $\eps$.
\end{proof}
\end{lem}
For its further use we recall that for every $\eps\in(0,1]$, the Malliavin derivative of the process $\{u^\eps(t,x), (t,x)\in[0,T]\times\Rd\}$ is an $\caH_T$-valued stochastic process $\{Du^\eps(t,x), (t,x)\in[0,T]\times\Rd\}$,
solution to
\begin{align}
\label{Dueps}
	Du^\eps(t,x)
	&	= \eps\Lambda(t-\cdot,x-\ast)\sigma(u^\eps(\cdot,\ast))\nonumber\\
 	&	+ \eps\int_0^t\int_{\Rd} \Lambda(t-s,x-y)\sigma^\prime(u^\eps(s,y))D u^\eps(s,y) M(\intd s,\intd y)\nonumber\\
	&	+ \int_0^t\int_{\Rd} \Lambda(t-s,x-y)b^\prime(u^\eps(s,y))D u^\eps(s,y) \intd y \intd s.
\end{align}
For the background on the Hilbert-valued stochastic and pathwise integrals in the preceding equation, we refer the reader to \cite{nualartquer} (see also \cite{dalangquer}).

\begin{lem}\label{lem:lem33}

Fix $(t,x)\in(0,T]\times\Rd$ and assume \ref{itm:finiteLambda}, \ref{itm:conditionic}, \ref{itm:boundsLambda}, $\sigma,b\in\caC^1$ with bounded derivatives and \ref{itm:sigmabound}. Then for every $p\in[1,\infty)$ there exists $C_{p}>0$ such that
\[ \big\|\gamma_{u^\eps(t,x)}^{-1}\big\|_{p} \leq C_{p}\eps^{-2}, \]
for any $\eps\in(0,1]$.
\begin{proof}
Fix $\eps\in(0,1]$ and $q\ge 2$. We will prove that there exists  $\zeta_0:=\zeta_0(q) >0$ such that for all $\zeta\in(0,\zeta_0)$
\begin{equation}
	P^\eps(\zeta) := \P\bigg\{\eps^{-2}\|Du^\eps(t,x)\|_{\caH_T}^2\leq \zeta\bigg\} \leq C\left(\zeta^q + \zeta^{\frac{2q\delta}{\gamma}}\right),
	\label{eq:prf331}
\end{equation}
where $C$ is a constant not depending on $\zeta$. 
Then, by the formula 
\begin{equation*}
\E(Y) = \int_0^\infty \P\{Y\geq \zeta\}\intd \zeta,
\end{equation*}
 valid for nonnegative random variables $Y$, the assertion will follow.

For any $0\le s<t$, let $\caH_{s,t}=L^2([s,t],\caH)$. Let $t_0$ be as defined in \ref{itm:boundsLambda}. We consider
 $\rho>0$ satisfying $\rho<t\wedge t_0$. From \eqref{Dueps} and the triangular inequality, we clearly have
\begin{align*}
  \|Du^\eps(t,x)\|_{\caH_T}^2
  & \geq \|Du^\eps(t,x)\|_{\caH_{t-\rho,t}}^2 \\
  & \geq \frac{1}{2}\|\eps\Lambda(t-\cdot,x-\ast)\sigma(u^\eps(\cdot,\ast))\|^2_{\caH_{t-\rho,t}} - \|X^\eps(t,x)\|_{\caH_{t-\rho,t}}^2,
\end{align*}
where
\begin{equation}
\label{dec}
X^\eps(t,x) := Du^\eps(t,x)-\eps\Lambda(t-\cdot,x-\ast)\sigma(u^\eps(\cdot,\ast)).
\end{equation}

The assumption \ref{itm:sigmabound} yields
\[ \|\eps\Lambda(t-\cdot,x-\ast)\sigma(u^\eps(\cdot,\ast))\|^2_{\caH_{t-\rho,t}} 
 \geq \eps^2\sigma_0^2 g_1(\rho). \]
Hence,
\begin{align}
	\P\bigg\{\eps^{-2}\|Du^\eps(t,x)\|_{\caH_T}^2\leq \zeta\bigg\}
	& \leq \P\bigg\{\eps^{-2}\|X^\eps(t,x)\|^2_{\caH_{t-\rho,t}} \geq \frac{\sigma_0^2}{2}g_1(\rho)-\zeta\bigg\} \notag \\
	& \leq \bigg(\frac{\sigma_0^2}{2}g_1(\rho)-\zeta\bigg)^{-q}\eps^{-2q}\E\big[\|X^\eps(t,x)\|_{\caH_{t-\rho,t}}^{2q}\big],
	\label{eq:proof341}
\end{align}
 where in the last inequality we have applied Chebyshev's inequality.

Our next objective is to find an upper bound for $\E\big[\|X^\eps(t,x)\|_{\caH_{t-\rho,t}}^{2q}\big]$. From \eqref{Dueps} we have
\begin{equation*}
\E\big[\|X^\eps(t,x)\|_{\caH_{t-\rho,t}}^{2q}\big]\le C\big(T_1(t,x; \rho,q)+T_2(t,x;\rho,q)\big),
\end{equation*}
with
\small
\begin{align*}
T_1^\eps(t,x; \rho,q)&=\E\left[\left\Vert \eps\int_0^t \int_{\Rd} \Lambda(t-s,x-y)\sigma^\prime(u^\eps(s,y)) Du^\eps(s,y) M(\intd s,\intd y)\right\Vert_{\caH_{t-\rho,t}}^{2q}\right],\\
T_2^\eps(t,x; \rho,q)&=\E\left[\left\Vert\int_0^t \int_{\Rd} \Lambda(t-s,x-y)b^\prime(u^\eps(s,y)) Du^\eps(s,y) \intd y \intd s\right\Vert_{\caH_{t-\rho,t}}^{2q}\right].
\end{align*}
\normalsize
The Malliavin derivative $D_{r,\ast}u^\eps(s,y)$ vanishes if $r\in(s,T]$. Thus, if $r\in[t-\rho,t]$, the domain of integration of the $s$ variable in the terms $T_1^\eps(t,x; \rho,q)$, $T_2^\eps(t,x; \rho,q)$ can be replaced by $[t-\rho,t]$. Moreover, following the proof of \cite[Lemma 8.2]{sanzbook}, we have
\begin{equation}
\label{3}
\sup_{\eps\in(0,1]}\sup_{s\in[0,\rho]}\sup_{y\in\Rd}\E\big[\|D_{t-\cdot,\ast} u^\eps(t-s,y)\|^{2q}_{\caH_{\rho}}\big] \le C(g_1(\rho))^q.
\end{equation}
By applying Burholder's inequality for Hilbert valued martingales (see for instance \cite{metivier}), we obtain
\begin{align}
\label{4}
	T_1^\eps(t,x; \rho,q)
	&	\le C\eps^{2q} (g_1(\rho))^q \sup_{\eps\in(0,1]}\sup_{s\in[t-\rho,t]}\sup_{y\in\Rd} \E\big[\|Du^\eps(s,y)\|^{2q}_{\caH_{t-\rho,t}}\big]	\nonumber\\
	&	\le C \eps^{2q} [g_1(\rho)]^{2q},
\end{align}
where in the last inequality we have used \eqref{3}.

We proceed now to the study of the term $T_2^\eps(t,x; \rho,q)$. For this, we use \eqref{dec} and Minkovski's inequality for the norm $\Vert\cdot\Vert_{\caH}$. So we are left with two terms that we study separately. For the first one, we use that $X_{r,\ast}(s,y)$ vanishes for $r\in(s,T]$, H\"older's inequality with respect to the finite measure $\Lambda(t-s,x-y)dsdy$, the boundedness of $b'$ and \eqref{eq:conditiondelta}. We obtain
\begin{align}
\label{5}
	&	\E\left[\left\Vert\int_0^t \int_{\Rd} \Lambda(t-s,x-y)b^\prime(u^\eps(s,y)) X^\eps(s,y) \intd y \intd s\right\Vert_{\caH_{t-\rho,t}}^{2q}\right]\nonumber\\
	&	\quad \le \E\left[\left(\int_{t-\rho}^t \int_{\Rd} \Lambda(t-s,x-y)|b^\prime(u^\eps(s,y))| \|X^\eps(s,y)\|_{\caH_{t-\rho,t}} \intd y \intd s\right)^{2q}\right]\nonumber \\
	& \quad \le C\rho^{\delta(2q-1)}\int_{t-\rho}^t \sup_{y\in\Rd} \E\big[\|X^\eps(s,y)\|_{\caH_{t-\rho,t}}^{2q}\big] J_2(t-s)ds \nonumber\\
	&	\quad \le C\rho^{\delta(2q-1)}\int_{t-\rho}^t \sup_{y\in\Rd} \E\big[\|X^\eps(s,y)\|_{\caH_{s-\rho,s}}^{2q}\big] J_2(t-s)ds,
\end{align}
where the last inequality follows from the property $\Vert \cdot\Vert_{\caH_{t-\rho,t}} \le \Vert \cdot\Vert_{\caH_{s-\rho,t}}$, for any $0\le s\le t\le T$.

Next we consider the second contribution from $T_2^\eps(t,x,\rho,q)$. The Lip\-schitz continuity of $\sigma$ together with Theorem \ref{basic} imply 
\small
\[ \sup_{{\eps\in(0,1]}}\sup_{(t,x)\in[0,T]\times\Rd} \E\big[|\sigma(u^\eps(t,x))|^q\big]
\leq C \bigg(1 + \sup_{{\eps\in(0,1]}}\sup_{(t,x)\in[0,T]\times\Rd} \E\big[|u^\eps(t,x)|^q\big]\bigg), \]
\normalsize
for any $q\in[1,\infty)$.
This yields
\begin{equation*}
	\sup_{y\in\Rd}\E\left[\Vert \Lambda(s-\cdot,y-\ast)\sigma(u^\eps(\cdot,\ast))\Vert^{2q}_{\caH_{s-\rho,s}}\right]\le C(g_1(\rho))^q.
\end{equation*}
Using this estimate and proceeding in a similar way as in the study of the previous term, we obtain
\begin{align}
\label{6}
	&	\E\left[\left\Vert\int_0^t \int_{\Rd} \Lambda(t-s,x-y)b^\prime(u^\eps(s,y))
\eps \Lambda(s-\cdot,y-\ast)\sigma(u^\eps(t,x)) \intd y \intd s\right\Vert_{\caH_{t-\rho,t}}^{2q}\right]\nonumber\\
	&	\qquad \le C \eps^{2q} \rho^{2q\delta} (g_1(\rho))^q.
\end{align}
With \eqref{4}, \eqref{5} and \eqref{6}, we have proved
\begin{align*}
	\sup_{x\in\Rd}\E\big[\|X^\eps(t,x)\|_{\caH_{t-\rho,t}}^{2q}\big]
	&	\le C \bigg(\eps^{2q}\left((g_1(\rho))^{2q} + \rho^{2q\delta}(g_1(\rho))^q\right)\\
	&	\quad + \rho^{\delta(2q-1)}\int_{t-\rho}^t  \sup_{y\in\Rd} 	\E\big[\|X^\eps(s,y)\|_{\caH_{s-\rho,s}}^{2q}\big]J_2(t-s) \intd s\bigg).
\end{align*}
Applying Gronwall's lemma in \cite[Lemma 15]{dalang} to the function 
\begin{equation*}
f(t)= \sup_{x\in\Rd}\E\big[\|X^\eps(t,x)\|_{\caH_{t-\rho,t}}^{2q}\big],
\end{equation*} we have
\begin{equation*}
	\sup_{x\in\Rd}\E\big[\|X^\eps(t,x)\|_{\caH_{t-\rho,t}}^{2q}\big]
	\le C \eps^{2q}\left((g_1(\rho))^{2q}+\rho^{2q\delta}(g_1(\rho))^q\right).
\end{equation*}
Plugging this estimate in \eqref{eq:proof341} we obtain
\[ \P\bigg\{\eps^{-2}\|Du^\eps(t,x)\|_{\caH_T}^2\leq \zeta\bigg\} \leq C\bigg(\frac{\sigma_0^2}{2}g_1(\rho)-\zeta\bigg)^{-q} \big((g_1(\rho))^{2q} + \rho^{2q\delta}(g_1(\rho))^q \big). \]
Let $0<\rho=\rho(\zeta)\le t\wedge t_0$ be such that $g_1(\rho)=\tfrac{4}{\sigma_0^2}\zeta$, which by \ref{itm:boundsLambda} implies $\rho\leq C\zeta^{1/\gamma}$. With this choice of $\rho$, the preceding inequality yields
\begin{equation*}
	\P\bigg\{\eps^{-2}\|Du^\eps(t,x)\|_{\caH_T}^2\leq \zeta\bigg\}
	\leq C\big(\zeta^q + \zeta^{\frac{2q\delta}{\gamma}}\big).
\end{equation*}

\end{proof}
\end{lem}
The proof of Theorem \ref{thm:LDPdensity} \ref{itm:LPDdensityi} is now complete.


\subsection{Lower bound}
\label{s3.2}
The purpose of this section is to prove that the family of random variables $(F^\eps)_{\eps\in(0,1]} = (u^\eps(t,x))_{\eps\in(0,1]}$, with fixed $(t,x)\in(0,T]\times \Rd$, satisfies the assumptions of  Proposition \ref{prop:bound441}, with $\psi(h)=\Phi_{t,x}^{h}$ (see \eqref{eq:defbhphi}), and we will identify the random variable $\caN$.  We will also prove that for any $h\in \caH_T$, $\bar{\gamma}_{\Phi_{t,x}^{h}}>0$.

\begin{lem}\label{lem:frechetdiff}
Assume \ref{itm:finiteLambda}, \ref{itm:conditionic} and that $\sigma,b\in\caC^1$ with Lipschitz continuous and bounded derivatives. Then, for all $(t,x)\in[0,T]\times\Rd$, the mapping
$\caH_T\ni h\mapsto \Phi^h_{t,x}$ defined in \eqref{eq:defbhphi} is Fr{\'e}chet differentiable.
\begin{proof}
Fix $(t,x)\in[0,T]\times\Rd$, $h\in\caH_T$. We use Cauchy-Schwarz' inequality, \ref{itm:conditionic}, \ref{itm:finiteLambda} and the Lipschitz continuity of $\sigma$ and $b$ to obtain
\begin{align*}
	|\Phi_{t,x}^{h}|^2
	& \leq C \left\{1 + \|h\|^2_{\caH_T} \|\Lambda(t-\cdot,x-\ast)\sigma(\Phi_{t,x}^{h})\|_{\caH_T}^2 \right.\\
	& \left. \qquad+ \bigg|\int_0^t\int_\Rd \Lambda(t-s,x-z) b(\Phi^{h}_{s,z})\intd z\intd s\bigg|^2\right\} \\
	& \leq C(\|h\|_{\caH_T}^2+1)\int_0^t \bigg(1+\sup_{(r,y)\in[0,s]\times\Rd} |\Phi^{h}_{r,y}|^2\bigg) \big(J_1(t-s)+J_2(t-s)\big)\intd s.
\end{align*}
Gronwall's Lemma yields
\begin{equation}
	\sup_{(t,x)\in[0,T]\times\Rd} |\Phi_{t,x}^{h}|^2 = C(\|h\|^2_{\caH_T}+1)\int_0^T\big(J_1(s)+J_2(s)\big)\intd s < \infty,
	\label{eq:finitebhphi}
\end{equation}
where the constant $C$ is independent of $h\in{\caH_T}$.
Now fix $h_0\in\caH_T$ and note that
\begin{align*}
  \Phi^{h+h_0}_{t,x} - \Phi_{t,x}^{h} = & \big\langle\Lambda(t-\cdot,x-\ast)\sigma(\Phi^{h+h_0}_{\cdot,\ast}),h_0\big\rangle_{\caH_T} \\
  & + \big\langle\Lambda(t-\cdot,x-\ast)\big(\sigma(\Phi^{h+h_0}_{\cdot,\ast})-\sigma(\Phi^{h}_{\cdot,\ast})\big),h\big\rangle_{\caH_T} \\
  & + \int_0^t\int_\Rd \Lambda(t-s,x-z)\big(b(\Phi^{h+h_0}_{s,z}) - b(\Phi^{h}_{s,z})\big)\intd z\intd s.
\end{align*}
With the same arguments as for the proof of \eqref{eq:finitebhphi}, we get that for all
$h_0\in \caH_T$
\begin{align*}
	& |\Phi^{h+h_0}_{t,x} - \Phi_{t,x}^{h}|^2 \leq C\|h_0\|_{\caH_T}^2\|\Lambda(t-\cdot,x-\ast)\sigma(\Phi^{h+h_0}_{\cdot,\ast})\|_{\caH_T}^2 \\
	& \qquad + C\big(\|h\|^2_{\caH_T}+1\big)\int_0^t \sup_{(r,y)\in[0,s]\times\Rd} |\Phi^{h+h_0}_{r,y} - \Phi^{h}_{r,y}|^2 \big(J_1(t-s)+J_2(t-s)\big)\intd s.
\end{align*}
Due to \eqref{eq:finitebhphi} and \ref{itm:finiteLambda}, the first term is bounded (up to a constant depending on $h$ and $h_0$) by $\|h_0\|^2_{\caH_T}$. Then applying Gronwall's Lemma we obtain
\begin{equation}
	\sup_{(t,x)\in[0,T]\times\Rd} |\Phi^{h+h_0}_{t,x} - \Phi_{t,x}^{h}| \leq C_{h,h_0}\|h_0\|_{\caH_T}.
	\label{eq:finitephi-phi}
\end{equation}
Note that the constant $C_{h,h_0}$ does not blow up as $\|h_0\|\to0$. With \eqref{eq:finitebhphi} and \eqref{eq:finitephi-phi}, we can prove the existence of the Fr{\'e}chet derivative of the map $h\mapsto\Phi_{t,x}^{h}$. First, we provide a candidate for it at the point ${g}\in \caH_T$, as follows:
\begin{align}
\label{9}
	\bar{D}\Phi_{t,x}^{h}({g})
	= & \big\langle \Xi^{h}(t,x),g\big\rangle_{\caH_T}\nonumber\\
	= & \big\langle \Lambda(t-\circ,x-\bullet)\sigma(\Phi^{h}_{\circ,\bullet}),g\big\rangle_{\caH_T} \nonumber\\
  & + \big\langle\Lambda(t-\cdot,x-\ast)\sigma'(\Phi^{h}_{\cdot,\ast})\big\langle\Xi^{h}(\cdot,\ast),g\big\rangle_{\caH_T},h\big\rangle_{\caH_T}\nonumber \\
  & + \int_0^t\int_\Rd \Lambda(t-s,x-z) b'(\Phi^{h}_{s,z})\big\langle\Xi^{h}(s,z),g\big\rangle_{\caH_T} \intd z\intd s,
\end{align}
where $\Xi^{h}(t,x)$ is defined by the  integral equation on ${\caH_T}$:
\begin{align}
	\Xi_{\circ,\bullet}^{h}(t,x) = & \Lambda(t-\circ,x-\bullet)\sigma(\Phi^{h}_{\circ,\bullet}) + \big\langle\Lambda(t-\cdot,x-\ast)\sigma'(\Phi^{h}_{\cdot,\ast})\Xi_{\circ,\bullet}^{h}(\cdot,\ast),h\big\rangle_{\caH_T} \notag\\
  & + \int_0^t\int_\Rd \Lambda(t-s,x-z) b'(\Phi^{h}_{s,z})\Xi_{\circ,\bullet}^{h}(s,z) \intd z\intd s. \label{eq:defXi}
\end{align}

According to Theorem \ref{thm:existenceanduniqueness}, this equation has a unique solution. Note that in the previous two formulas $(\circ,\bullet)$ is the argument in $[0,T]\times\Rd$ which interacts with ${g}\in \caH_T$ (the element at which $\Xi^{h}(t,x)$ is evaluated) and $(\cdot,\ast)$ is the argument in $[0,T]\times\Rd$ that interacts with $h\in{\caH_T}$  which is the point where the Fr{\'e}chet derivative is taken.

From \eqref{9} and \eqref{eq:defXi}, we clearly have
\begin{align}
& \frac{\Phi^{h+h_0}_{t,x} - \Phi_{t,x}^{h} - \bar{D}\Phi_{t,x}^{h}(h_0)}{\|h_0\|_{\caH_T}} \notag\\
 = & \frac{1}{\|h_0\|_{\caH_T}} \big\langle\Lambda(t-\circ,x-\bullet)\big(\sigma(\Phi^{h+h_0}_{\circ,\bullet}) - \sigma(\Phi^{h}_{\circ,\bullet})\big),h_0\big\rangle_{\caH_T} \notag\\
	& + \bigg\langle\Lambda(t-\cdot,x-\ast)\bigg(\frac{\sigma(\Phi^{h+h_0}_{\cdot,\ast}) - \sigma(\Phi^{h}_{\cdot,\ast}) - \sigma'(\Phi^{h}_{\cdot,\ast})\bar{D}\Phi^{h}_{\cdot,\ast}(h_0)}{\|h_0\|_{\caH_T}}\bigg),h\bigg\rangle_{\caH_T} \notag\\
	& + \int_0^t\int_\Rd \Lambda(t-s,x-z)\frac{b(\Phi^{h+h_0}_{s,z}) - b(\Phi^{h}_{s,z}) - b'(\Phi^{h}_{s,z})\bar{D}\Phi_{t,x}^{h}(h_0)}{\|h_0\|_{\caH_T}} \intd z\intd s. \label{eq:frechetder}
\end{align}

Our aim is to have an upper bound for the absolute value of each term on the right-hand side of \eqref{eq:frechetder}.
By applying  Cauchy-Schwarz' inequality, the fact that $\sigma$ is Lipschitz continuous, \ref{itm:finiteLambda} and \eqref{eq:finitephi-phi}, we have
\begin{align*}
\frac{1}{\|h_0\|^2_{\caH_T}} & \big\vert\big\langle\Lambda(t-\circ,x-\bullet)\big(\sigma(\Phi^{h+h_0}_{\circ,\bullet}) - \sigma(\Phi^{h}_{\circ,\bullet})\big),h_0\big\rangle_{\caH_T}^2\big\vert \\
& \leq \big\|\Lambda(t-\cdot,x-\ast)\big(\sigma(\Phi^{h+h_0}_{\cdot,\ast}) - \sigma(\Phi^{h}_{\cdot,\ast})\big)\big\|_{\caH_T}^2 \\
& \leq \sup_{(r,y)\in[0,T]\times\Rd} \big|\sigma(\Phi^{h+h_0}_{r,y}) - \sigma(\Phi^{h}_{r,y})\big|^2 \|\Lambda(t-\cdot,x-\ast)\|^2_{\caH_T} \\
& \leq C_{h,h_0}\|h_0\|^2_{\caH_T}.
\end{align*}
For the second term, we first use Cauchy-Schwarz' inequality and apply the usual procedure involving the Fourier transformation. Then we use the mean-value theorem to obtain
\small
\begin{align*}
	& \Bigg|\bigg\langle\Lambda(t-\cdot,x-\ast)\bigg(\frac{\sigma(\Phi^{h+h_0}_{\cdot,\ast}) - \sigma(\Phi^{h}_{\cdot,\ast}) - \sigma'(\Phi^{h}_{\cdot,\ast})\bar{D}\Phi^{h}_{\cdot,\ast}(h_0)}{\|h_0\|_{\caH_T}}\bigg),h\bigg\rangle_{\caH_T} \Bigg|^2 \notag\\
	\leq & \|h\|^2_{\caH_T}\int_0^t \sup_{(r,y)\in[0,s]\times\Rd} \bigg|\frac{\sigma(\Phi^{h+h_0}_{r,y}) - \sigma(\Phi^{h}_{r,y}) - \sigma'(\Phi^{h}_{r,y})\bar{D}\Phi^{h}_{r,y}(h_0)}{\|h_0\|_{\caH_T}}\bigg|^2 J_1(t-s)\intd s \notag\\
	= & \|h\|^2_{\caH_T}\int_0^t \sup_{(r,y)\in[0,s]\times\Rd} \bigg|\frac{\sigma'(\xi^{h,h_0}_{r,y})(\Phi^{h+h_0}_{r,y} - \Phi^{h}_{r,y}) - \sigma'(\Phi^{h}_{r,y})\bar{D}\Phi^{h}_{r,y}(h_0)}{\|h_0\|_{\caH_T}}\bigg|^2 J_1(t-s)\intd s \notag\\
	\leq & \|h\|^2_{\caH_T}\int_0^t \sup_{(r,y)\in[0,s]\times\Rd} \bigg|\big(\sigma'(\xi^{h,h_0}_{r,y})-\sigma'(\Phi^{h}_{r,y})\big)\frac{\Phi^{h+h_0}_{r,y} - \Phi^{h}_{r,y}}{\|h_0\|_{\caH_T}}\bigg|^2 J_1(t-s)\intd s \notag\\
	& + \|h\|^2_{\caH_T}\int_0^t\sup_{(r,y)\in[0,s]\times\Rd} \bigg|\sigma'(\Phi^{h}_{r,y})\frac{\Phi^{h+h_0}_{r,y} - \Phi^{h}_{r,y} - \bar{D}\Phi^{h}_{r,y}(h_0)}{\|h_0\|_{\caH_T}}\bigg|^2 J_1(t-s)\intd s,
\end{align*}
\normalsize
where $\xi^{h,h_0}_{r,y}$ is a real number in the convex hull of $\Phi^{h+h_0}_{r,y}$ and $\Phi^{h}_{r,y}$.

Using \eqref{eq:finitephi-phi} and the Lipschitz continuity property of $\sigma'$, along with \ref{itm:boundsLambda},
the first term on the right-hand side of the last inequality can be bounded from above by
\[ C\|h\|_{\caH_T}^2 \sup_{(t,x)\in[0,T]\times\Rd} \big|\xi^{h,h_0}_{t,x} - \Phi_{t,x}^{h}\big|^2, \]
and therefore also by $C_{h,h_0}\|h_0\|_{\caH_T}^2\|h\|_{\caH_T}^2$.

We are assuming that $\sigma'$ is bounded. Hence, we have proved
\begin{align*}
& \Bigg|\bigg\langle\Lambda(t-\cdot,x-\ast)\bigg(\frac{\sigma(\Phi^{h+h_0}_{\cdot,\ast}) - \sigma(\Phi^{h}_{\cdot,\ast}) - \sigma'(\Phi^{h}_{\cdot,\ast})\bar{D}\Phi^{h}_{\cdot,\ast}(h_0)}{\|h_0\|_{\caH_T}}\bigg),h\bigg\rangle_{\caH_T} \Bigg|^2 \notag\\
& \leq C_{h,h_0}\|h_0\|_{\caH_T} + C_h\int_0^t \sup_{(r,y)\in[0,s]\times\Rd}\bigg|\frac{\Phi^{h+h_0}_{r,y} - \Phi^{h}_{r,y} - \bar{D}\Phi^{h}_{r,y}(h_0)}{\|h_0\|_{\caH_T}}\bigg|^2 J_1(t-s)\intd s.
\end{align*}
A similar estimate, with  $J_1(t-s)$ replaced by $J_2(t-s)$ holds for the last term on the right-hand side of \eqref{eq:frechetder}.

Summarizing, we have proved that
\begin{align*}
	& \bigg|\frac{\Phi^{h+h_0}_{t,x} - \Phi_{t,x}^{h} - \bar{D}\Phi_{t,x}^{h}(h_0)}{\|h_0\|_{\caH_T}}\bigg|^2 \\
	& \leq C_{h,h_0}\|h_0\|^2_{\caH_T} + C_h\int_0^t \sup_{(r,y)\in[0,s]\times\Rd}\bigg|\frac{\Phi^{h+h_0}_{r,y} - \Phi^{h}_{r,y} - \bar{D}\Phi^{h}_{r,y}(h_0)}{\|h_0\|_{\caH_T}}\bigg|^2\\
	 &\qquad \times \big(J_1(t-s) + J_2(t-s)\big)\intd s.
\end{align*}
Since for any $h\in\caH_T$, $\sup_{\|h_0\|_{\caH_T}\le 1}C_{h,h_0} <\infty$,
by Gronwall's Lemma we conclude
\begin{equation*}
\lim_{\|h_0\|_{\caH_T}\to 0}\bigg|\frac{\Phi^{h+h_0}_{t,x} - \Phi_{t,x}^{h} - \bar{D}\Phi_{t,x}^{h}(h_0)}{\|h_0\|_{\caH_T}}\bigg|=0.
\end{equation*}
This ends the proof of the Lemma.
\end{proof}
\end{lem}

\begin{numrem}
Assume \ref{itm:sigmabCinf}. By further differentiating the term $\Xi^{h}(t,x)$ in \eqref{eq:defXi} and repeating the calculation involving the definition of higher-order Fr{\'e}\-chet differentiability, it can be shown that $\Phi_{t,x}$ is Fr{\'e}chet differentiable of any order.
\end{numrem}

Now we are in position to check \eqref{eq:prop441}.

\begin{lem}\label{lem:lem45}
Fix ${(t,x)\in[0,T]\times\Rd}$ and assume \ref{itm:finiteLambda}, \ref{itm:conditionic}, \ref{itm:boundsLambda}, $\sigma,b\in\caC^1$ with Lipschitz continuous and bounded derivatives, and \ref{itm:sigmabound}. Then, for all $h\in \caH_T$, \eqref{eq:prop441} holds with $F^\eps = u^\eps(t,x)$, $\psi(h) = \Phi_{t,x}^{h}$ and
 $\caN_{t,x}(h)$ given by the  SPDE
\begin{align}
\label{n}
  \caN_{t,x}(h) = & \int_0^t\int_\Rd \Lambda(t-s,x-y)\sigma(\Phi^{h}_{s,y})M(\intd s,\intd y)\nonumber\\
  & + \langle\Lambda(t-\cdot,x-\ast)\sigma'(\Phi^{h}_{\cdot,\ast})
  \caN_{\cdot,\ast}(h),h\rangle_{\caH_T}\nonumber \\
  & + \int_0^t\int_\Rd \Lambda(t-s,x-y) b'(\Phi^{h}_{s,y})\caN_{s,y}(h) \intd y\intd s.
\end{align}
\begin{proof}
First we note that by Theorem \ref{thm:existenceanduniqueness} there exists a unique solution to \eqref{n}. 
The integrand in the stochastic integral term of this equation is deterministic, consequently, the random variable $\caN_{t,x}(h)$ is Gaussian and therefore it belongs to the first Wiener chaos. Its variance is given by $\|D\caN_{t,x}(h)\|^2_{\caH_T}$, where $D$ denotes the Malliavin derivative.

The Malliavin derivative of $\caN_{t,x}(h)$ satisfies the equation
\begin{align}
\label{eq:DN}
  D\caN_{t,x}(h) = & \Lambda(t-\cdot,x-\ast)\sigma(\Phi^{h}_{\cdot,\ast})
  + \langle \Lambda(t-\cdot,x-\ast)\sigma'(\Phi^{h}_{\cdot,\ast})D\caN_{\cdot,\ast}(h) , h \rangle_{\caH_T}\nonumber \\
  & + \int_0^t\int_\Rd \Lambda(t-s,x-y) b'(\Phi^{h}_{s,y})D\caN_{s,y}(h) \intd y\intd s.
\end{align}
Comparing this equation with the one for $\bar{D}\Phi^h_{t,x}$ in \eqref{9}, \eqref{eq:defXi} and invoking the uniqueness of solution, we see that, for any $h\in\caH_T$, the $\caH_T$-valued stochastic processes
$\{D\caN_{t,x}(h), (t,x)\in[0,T]\times \R^d\}$ and $\{\bar{D}\Phi^h_{t,x}, (t,x)\in[0,T]\times \R^d\}$ are indistinguishable.
In particular, the variance of $\caN_{t,x}(h)$ is $\|\bar{D}\Phi^h_{t,x}\|_{\caH_T}^2$.

Set $u^{\eps,h}(t,x):=u(t,x;\omega+\eps^{-1}h)$. According to Lemma \ref{cor:defuh}, the process $(u^{\eps,h}(t,x), (t,x)\in[0,T]\times \Rd)$ satisfies \eqref{eq:defuh}. By uniqueness of solution we clearly have $u^{\eps,0}(t,x) = u^\eps(t,x)$ and $u^{0,h}(t,x) = \lim_{\eps\downarrow0} u^{\eps,h}(t,x) = \Phi_{t,x}^{h}$, for any $(t,x)\in[0,T]\times \Rd$.

Next, we prove in our context the convergence \eqref{eq:prop441} in $L^p(\Omega)$ norm, for $p\in[2,\infty)$. Set
\begin{equation*}
Z^{\eps,h}_{t,x} = \eps^{-1}(u^{\eps,h}(t,x)-\Phi_{t,x}^{h})-\caN_{t,x}(h).
\end{equation*}
By using the equations satisfied by each one of the terms on the right hand-side of that expression, we see that
\begin{equation*}
\E\big[|Z^{\eps,h}_{t,x}|^p\big]\le C\sum_{i=1}^3 T^{\eps,h,i}_{t,x},
\end{equation*}
where
\small
\begin{align*}
T^{\eps,h,1}_{t,x}  & = \E\Bigg[\bigg|\int_0^t\int_\Rd \Lambda(t-s,x-z)\big(\sigma(u^{\eps,h}(s,z))-\sigma(\Phi^{h}_{s,z})\big)M(\intd s,\intd z)\bigg|^p\Bigg], \notag\\
T^{\eps,h,2}_{t,x} & = \E\Bigg[\bigg|\bigg\langle \Lambda(t-\cdot,x-\ast)\bigg(\frac{\sigma(u^{\eps,h}(\cdot,\ast)) - \sigma(\Phi^{h}_{\cdot,\ast})}{\eps} - \sigma'(\Phi^{h}_{\cdot,\ast})\caN_{\cdot,\ast}(h)\bigg),h\bigg\rangle_{\caH_T}\bigg|^p\Bigg], \notag\\
T^{\eps,h,3}_{t,x}  & = \E\Bigg[\bigg|\int_0^t\int_\Rd \Lambda(t-s,x-z)\bigg(\frac{b(u^{\eps,h}(s,z)) - b(\Phi^{h}_{s,z})}{\eps} - b'(\Phi^{h}_{s,z})\caN_{s,z}(h)\bigg)\intd z\intd s\bigg|^p\Bigg].
\end{align*}
\normalsize
We will prove that each one of these terms tends to zero as $\eps\downarrow0$.

By the usual estimates on moments of stochastic and pathwise integrals, we have
\begin{align*}
	&\sup_{(r,y)\in[0,t]\times\Rd} \E\big[|u^{\eps,h}(r,y)-\Phi^{h}_{r,y}|^p\big] \notag\\
	& \qquad\leq C\eps^p\bigg(1+\sup_{\eps\in(0,1]}\sup_{(t,x)\in[0,T]\times\Rd} \E\big[|u^{\eps,h}(t,x)|^p\big]\bigg) \notag\\
	& \qquad\phantom{\leq} + C\int_0^t \sup_{(r,y)\in[0,s]\times\Rd} \E\big[|u^{\eps,h}(r,y)-\Phi^{h}_{r,y}|^p\big] \big(J_1(t-s)+J_2(t-s)\big)\intd s.
\end{align*}
By Gronwall's Lemma this yields
\begin{equation}
\label{eq:Lplimitofuepsh}
\lim_{\eps\downarrow0}\left(\sup_{(t,x)\in[0,T]\times\Rd} \E\big[|u^{\eps,h}(t,x) - \Phi_{t,x}^{h}|^p\big]\right)=0.
\end{equation}
Since
\begin{equation*}
T^{\eps,h,1}_{t,x} \le C \sup_{(t,x)\in[0,T]\times\Rd} \E\big[|u^{\eps,h}(t,x) - \Phi_{t,x}^{h}|^p\big],
\end{equation*}
we deduce
\begin{equation}
\label{10}
\lim_{\eps\downarrow0}\left(\sup_{(t,x)\in[0,T]\times\Rd} T^{\eps,h,1}_{t,x}\right)=0.
\end{equation}
Next, we deal with the term $T^{\eps,h,2}_{t,x}$. Cauchy-Schwarz's inequality and the mean-value theorem applied to $\sigma$ yield
\small
\begin{align}
	& \E\Bigg[\bigg|\bigg\langle\Lambda(t-\cdot,x-\ast)\bigg(\frac{\sigma(u^{\eps,h}(\cdot,\ast)) - \sigma(\Phi^{h}_{\cdot,\ast})}{\eps} - \sigma'(\Phi^{h}_{\cdot,\ast})\caN_{\cdot,\ast}(h)\bigg),h\bigg\rangle_{\caH_T}\bigg|^p\Bigg] \notag \\
	& \leq C\|h\|^p_{\caH_T}\nonumber\\
	&\quad \times\int_0^t \sup_{(r,y)\in[0,s]\times\Rd} \E\Bigg[\bigg|\frac{\sigma(u^{\eps,h}(r,y)) - \sigma(\Phi^{h}_{r,y})}{\eps} - \sigma'(\Phi^{h}_{r,y})\caN_{r,y}(h)\bigg|^p\Bigg] J_1(t-s)\intd s \notag\\
	& \leq C\left\{\int_0^t \sup_{(r,y)\in[0,s]\times\Rd} \E\Bigg[\bigg|\sigma'(\xi^{\eps,h}_{r,y})\bigg(\frac{u^{\eps,h}(r,y) - \Phi^{h}_{r,y}}{\eps} - \caN_{r,y}(h)\bigg)\bigg|^p\Bigg] J_1(t-s)\intd s \notag\right.\\
	& \left.\phantom{\leq} + \sup_{(t,x)\in[0,T]\times\Rd} \E\Big[\big|\big(\sigma'(\Phi_{t,x}^{h}) - \sigma'(\xi^{\eps,h}_{t,x})\big)\caN_{t,x}(h)\big|^p\Big] \int_0^t J_1(t-s)\intd s\right\}, \label{eq:prf451}
\end{align}
\normalsize
where $\xi^{\eps,h}_{r,y}(\omega)$ is a point lying in the open interval determined by  $\Phi^{h}_{r,y}$ and $u^{\eps,h}(r,y;\omega)$.

From Theorem \ref{thm:existenceanduniqueness} and the Lipschitz continuity of $\sigma^\prime$ we have
\begin{align*}
& \sup_{(t,x)\in[0,T]\times\Rd} \E\Big[\big|\big(\sigma'(\Phi_{t,x}^{h}) - \sigma'(\xi^{\eps,h}_{t,x})\big)\caN_{t,x}(h)\big|^p\Big] \notag\\
            & \leq \sup_{(t,x)\in[0,T]\times\Rd}\left( \E\Big[\big|\sigma'(\Phi_{t,x}^{h}) - \sigma'(\xi^{\eps,h}_{t,x})\big|^{2p}\Big]\right)^{\frac12} \sup_{(t,x)\in[0,T]\times\Rd} \left(\E\Big[\big|\caN_{t,x}(h)\big|^{2p}\Big]\right)^{\frac12} \notag\\
            & \leq C\sup_{(t,x)\in[0,T]\times\Rd} \left(\E\big[\big|\Phi_{t,x}^{h} - u^{\eps,h}_{t,x}\big|^{2p}\big]\right)^{1/2}.
\end{align*}
Consequently,
\begin{align}
\label{11}
T^{\eps,h,2}_{t,x}&\le C\sup_{(t,x)\in[0,T]\times\Rd} \left(\E\big[\big|\Phi_{t,x}^{h} - u^{\eps,h}_{t,x}\big|^{2p}\big]\right)^{1/2}\nonumber\\
& + C\int_0^t \sup_{(r,y)\in[0,s]\times\Rd} \E\Bigg[\bigg|\frac{u^{\eps,h}(r,y) - \Phi^{h}_{r,y}}{\eps} - \caN_{r,y}(h)\bigg|^p\Bigg] J_1(t-s)\intd s,
\end{align}
 With similar arguments, one can check that
 \begin{align}
 \label{110}
T^{\eps,h,3}_{t,x}\le & C\sup_{(t,x)\in[0,T]\times\Rd} \left(\E\big[\big|\Phi_{t,x}^{h} - u^{\eps,h}_{t,x}\big|^{2p}\big]\right)^{1/2}\nonumber\\
& + C\int_0^t \sup_{(r,y)\in[0,s]\times\Rd} \E\Bigg[\bigg|\frac{u^{\eps,h}(r,y) - \Phi^{h}_{r,y}}{\eps} - \caN_{r,y}(h)\bigg|^p\Bigg] J_2(t-s)\intd s,
\end{align}

Thus from \eqref{10}, \eqref{11}, \eqref{110} it follows that
\begin{align*}
	\sup_{(r,y)\in[0,t]\times\Rd} & \E\big[|Z^{\eps,h}_{r,y}|^p\big] \\
	& \leq C_\eps + C\int_0^t \sup_{(r,y)\in[0,s]\times\Rd} \E\big[|Z^{\eps,h}_{r,y}|^p\big] \big(J_1(t-s)+J_2(t-s)\big)\intd s,
 \end{align*}
where  $C_\eps$ converges to zero as $\eps\downarrow0$. Applying Gronwall's Lemma we see that $Z^{\eps,h}_{t,x}$ converges to zero in $L^p$ as $\eps\downarrow0$ for all $h\in \caH_T$, uniformly in ${(t,x)\in[0,T]\times\Rd}$.

The next step consists of proving the convergence to zero of  $Z^{\eps,h}_{t,x}$ in the $\D^{1,p}$ norm, for any $p\in[2,\infty)$. Since $\Phi_{t,x}^{h}$ is deterministic this reduces to show that $\eps^{-1}Du^{\eps}(t,x;\omega+h) - D\caN(h)$ converges to zero as $\eps\downarrow 0$ in $L^p(\Omega;{\caH_T})$.

By applying the Malliavin derivative operator to Equation \eqref{eq:defuh}, one can show that the process $Du^{\eps,h}(t,x)= Du^{\eps}(t,x;\omega+h)$ satisfies the SPDE
\begin{align}
\label{111}
  Du^{\eps,h}(t,x)
  = & \eps\Lambda(t-\cdot,x-\ast)\sigma(u^{\eps,h}(\cdot,\ast))\nonumber \\
  & + \eps\int_0^t\int_\Rd \Lambda(t-s,x-z)\sigma'(u^{\eps,h}(s,z))Du^{\eps,h}(s,z)M(\intd s,\intd z)\nonumber \\
  & + \langle\Lambda(t-\cdot,x-\ast)\sigma'(u^{\eps,h}(\cdot,\ast))Du^{\eps,h}(\cdot,\ast),h\rangle_{\caH_T}\nonumber \\
  & + \int_0^t\int_\Rd \Lambda(t-s,x-z) b'(u^{\eps,h}(s,z))Du^{\eps,h}(s,z) \intd z\intd s.
\end{align}
For its further use, we remark that
\begin{equation}
\label{13}
\lim_{\eps\downarrow0}\left(\sup_{(t,x)\in[0,T]\times \Rd}\E\big[\|Du^{\eps,h}(t,x)\|_{\caH_T}^p\big]\right)=0.
\end{equation}
Indeed, this follows from the estimate
\begin{align*}
	& \E\big[\|Du^{\eps,h}(t,x)\|_{\caH_T}^p\big]
\le  C\eps^p\bigg(1+\sup_{{\eps\in(0,1]}}\sup_{(t,x)\in[0,T]\times\Rd} \E\big[|u^{\eps,h}(t,x)|^p\big]\bigg)\nonumber\\
 &\qquad + C\int_0^t \sup_{{\eps\in(0,1]}}\sup_{(r,y)\in[0,s]\times\Rd} \E\big[\|Du^{\eps,h}(r,y)\|_{\caH_T}^p\big] \big(J_1(t-s)+J_2(t-s)\big)\intd s,
 \end{align*}
along with Gronwall's lemma.

By \eqref{eq:DN} and \eqref{111} we easily obtain
\small
\begin{align}
\label{12}
  & \E\big[\|\eps^{-1}Du^{\eps,h}(t,x) - D\caN_{t,x}(h)\|_{\caH_T}^p\big] \nonumber\\
  & \leq C\E\big[\big\|\Lambda(t-\cdot,x-\ast)\big(\sigma(u^{\eps,h}(\cdot,\ast))-\sigma(\Phi^{h}_{\cdot,\ast})\big)\big\|_{\caH_T}^p\big] \nonumber\\
  & + C\E\Bigg[\bigg\|\int_0^t\int_\Rd \Lambda(t-s,x-z) \sigma'(u^{\eps,h}(s,z))Du^{\eps,h}(s,z) M(\intd s,\intd z)\bigg\|_{\caH_T}^p\Bigg] \nonumber\\
  & + C\E\Bigg[\bigg\|\big\langle\Lambda(t-\cdot,x-\ast)\big(\eps^{-1}\sigma'(u^{\eps,h}(\cdot,\ast))Du^{\eps,h}(\cdot,\ast) - \sigma'(\Phi^{h}_{\cdot,\ast})D\caN_{\cdot,\ast}(h)\big),h\big\rangle_{\caH_T}\bigg\|_{\caH_T}^p\Bigg] \nonumber\\
  & + C\E\Bigg[\bigg\|\int_0^t\int_\Rd \Lambda(t-s,x-z)\big(\eps^{-1}b'(u^{\eps,h}(s,z))Du^{\eps,h}(s,z) - b'(\Phi^{h}_{s,z})\nonumber \\
  & \phantom{mmmmmm} \times D\caN_{s,z}(h)\big)\intd z\intd s\bigg\|_{\caH_T}^p\Bigg].
\end{align}
\normalsize
Each term on the right-hand side of \eqref{12} converges to zero as $\eps\downarrow 0$.
Indeed, for the first and second terms, this is a consequence of \eqref{eq:Lplimitofuepsh} and \eqref{13}, respectively.
For the analysis of the last two ones, we use the argument involving the mean-value Theorem as in \eqref{eq:prf451}.
Then Gronwall's Lemma yields the assertion.

In order to finish the proof, we must check the convergence of $Z^{\eps,h}_{t,x}$ in the $\D^{k,p}$ norm, for any $k\ge 2$, $p\in[2,\infty)$.
Since $\caN_{t,x}(\bar h)$ is a Gaussian random variable, this reduces to show that
$\eps^{-1}D^ku^{\eps,h}(t,x)$ converges to zero in $L^p(\Omega;{\caH_T}^{\otimes k})$, which is proved recursively on $k\ge 2$. We leave the details to the reader.
\end{proof}
\end{lem}

Thanks to Proposition \ref{prop:bound441}, the results proved so far establish the lower bound in \eqref{eq:ldpdenslower} with
\begin{equation*}
 I(y) = \inf \bigg\{\frac{1}{2}\|h\|_{\caH_T}^2; \; h\in \caH_T, \Phi_{t,x}^{h} = y, \bar{\gamma}_{\Phi_{t,x}^{h}}>0 \bigg\}.
\end{equation*}
In the next lemma it is shown that under the standing assumptions, the condition $\bar\gamma_{\psi}(h)>0$ is satisfied. Hence, $I(y)$ is as in \eqref{eq:ratefunc}.

\begin{lem}\label{lem:44}
Fix $(t,x)\in(0,T]\times\Rd$ and assume \ref{itm:finiteLambda}, \ref{itm:conditionic}, \ref{itm:boundsLambda}, $\sigma,b\in\caC^1$ with bounded derivatives and \ref{itm:sigmabound}. Then $\bar{\gamma}_{\Phi_{t,x}^{h}}>0$ for all $h\in \caH_T$.
\begin{proof}
The proof follows the same strategy as in Lemma \ref{lem:lem33}, the difference being that here we use deterministic arguments.

 Fix $h\in \caH_T$ and let $0<\rho<t\wedge t_0$, with $t_0$ as in \ref{itm:boundsLambda}. Remember that  $\bar{\gamma}_{\Phi_{t,x}^{h}} = \Vert \bar{D}\Phi_{t,x}^{h}\Vert_{\caH_T}^2$, where $\bar{D}$ stands for the Fr\'echet derivative. Using \eqref{9} and \eqref{eq:defXi} we clearly obtain
 \begin{equation*}
 \bar{\gamma}_{\Phi_{t,x}^{h}} = \|\Xi^{h}_{\circ,\bullet}(t,x)\|^2_{\caH_T} \geq
 \|\Xi^{h}_{r,\bullet}(t,x)\|^2_{\caH_{t-\rho,t}}
 \geq \frac{1}{2}A^1_{t,x}(\rho)-A^2_{t,x}(\rho),
 \end{equation*}
 with
 \begin{align*}
	A^1_{t,x}(\rho) & = \|\Lambda(t-\circ,\bullet)\sigma(\Phi^{h}_{\circ,\bullet})\|^2_{\caH_{t-\rho,t}},\\
	A^2_{t,x}(\rho) & = \|\chi_{\circ,\bullet}(t,x)\|_{\caH_{t-\rho,t}}^2
\end{align*}
and
\[ \chi_{\circ,\bullet}(t,x) = \Xi^{h}_{\circ,\bullet}(t,x) - \Lambda(t-\circ,\bullet)\sigma(\Phi^{h}_{\circ,\bullet}). \]
Using \ref{itm:boundsLambda}, we have
\begin{equation}
	A^1_{t,x}(\rho) \geq \sigma_0^2g_1(\rho),
	\label{eq:prflem381}
\end{equation}
and this bound is uniform in ${(t,x)\in[0,T]\times\Rd}$ and in $h\in \caH_T$.

Our next aim is to prove that there exists $\zeta>0$ and $0<\rho<t\wedge t_0$ such that $\tfrac{\sigma^2_0}{2}g_1(\rho) - A^2_{t,x}(\rho) \geq \zeta$, or equivalently,
\begin{equation}
	\bigg(\frac{\sigma^2_0}{2}g_1(\rho)-\zeta\bigg)^{-1}A^2_{t,x}(\rho) \leq 1.
	\label{eq:proof441}
\end{equation}
For this, we will find a suitable upper bound for $A^2_{t,x}(\rho)$. By using the definition of $\chi_{\circ,\bullet}(t,x)$ and \eqref{eq:defXi}, we have
\begin{equation*}
 A^2_{t,x}(\rho) \leq C \left(N^1_{t,x}(\rho)+N^2_{t,x}(\rho)\right),
 \end{equation*}
 with
 \begin{align*}
 & N^1_{t,x}(\rho)= \left\Vert\big\langle\Lambda(t-\cdot,x-\ast)\sigma'(\Phi^{h}_{\cdot,\ast})\Xi^{h}(\cdot,\ast),h\big\rangle_{\caH_T} \right\Vert^2_{\caH_{t-\rho,t}}, \nonumber\\
 & N^2_{t,x}(\rho)= \left\Vert\int_0^t\int_\Rd \Lambda(t-s,x-z) b'(\Phi^{h}_{s,z})\Xi^{h}(s,z) \intd z\intd s\right\Vert_{\caH_{t-\rho,t}}^2.
 \end{align*}
 Remember that, similarly as in \eqref{eq:defXi}, $(\circ,\bullet)$ is the argument in $[t-\rho,t]\times\Rd$ relevant for the
 $\caH_{t-\rho,t}$ norm, while $(\cdot,\ast)$ interacts with $h$ in the $\caH_T$ norm.

From the definition of $\Xi(s,y)$ given in \eqref{eq:defXi} and Proposition \ref{A.3}, we deduce the estimate
\[ \sup_{y\in\R^d} \| \Xi_{\circ,\bullet}^{h}(s,y)\|_{\caH_{s-\rho,s}}^2 \leq C g_1(\rho). \]
By applying first Schwarz's inequality to the inner product in $\caH_T$, the preceding estimate yields
\begin{equation}
\label{20}
	N^1_{t,x}(\rho)\le C (g_1(\rho))^2.
\end{equation}
(See \eqref{4} for an analogous result).

By similar arguments, we have
\begin{equation}
\label {21}
N^2_{t,x}(\rho)\le C \rho^{2\delta} g_1(\rho).
\end{equation}
(Notice the analogy with \eqref{5}).

With \eqref{20}, \eqref{21}, we see that the left-hand side of \eqref{eq:proof441} is bounded by
\begin{equation*}
	C\bigg(\frac{\sigma^2_0}{2}g_1(\rho)-\zeta\bigg)^{-1}\big(g_1(\rho)^2 + \rho^{2\delta} g_1(\rho)\big).
\end{equation*}
Fix $\rho=\rho(\zeta)>0$ such that $g_1(\rho)=\tfrac{4}{\sigma_0^2}\zeta$ which by \ref{itm:boundsLambda} implies that $\rho\leq C\zeta^{1/\gamma}$. Then the previous expression is bounded by $C(\zeta+\zeta^{\frac{2\delta}{\gamma}})$. Hence, \eqref{eq:proof441} holds
for a suitable choice of $\zeta>0$. This ends the proof.
\end{proof}
\end{lem}

The results established so far show that Theorem \ref{thm:LDPdensity} (ii) holds for the set $y\in\R$ such that $p^\eps_{t,x}(y)>0$ for all $\eps$ small enough. The next objective is to analyze when this condition is satisfied and also whether the function $I$
defined in \eqref{eq:ldpdenslower} is finite.  Both questions are related to the characterisation of the topological support of the law of the random variable $u^{\eps}(t,x)$, with $\eps\in(0,1]$ and $(t,x)\in(0,T]\times \R^d$ fixed.

Under suitable conditions, we prove in Theorem \ref{ta.1} that the support of $u^\eps(t,x)$ does not depend on the parameter $\eps$ and is given by
\[ \mathcal{S}:=\supp(P\circ[u^\eps(t,x)]^{-1}) = \overline{\{\Phi_{t,x}^{h};\; h\in \caH_T\}}. \]
In particular, $\mathcal{S}$ is the topological support of the law of the random variable $u^1(t,x)$, that we denote by $u(t,x)$.

Since for any $p\in[2,\infty)$, $u^\eps(t,x)\in \mathbb{D}^{1,p}$ (see Lemma \ref{lem:nondegenerateI}), we can apply Fang's result quoted in  \cite[Proposition 4.1.1]{nualartstflour} to deduce that $\mathcal{S}$ is a closed interval. Moreover, applying \cite[Proposition 4.1.2]{nualartstflour}, we obtain that for all points in the interior of $\mathcal{S}$, denoted by $\mathring{\mathcal{S}}$, we have $p^\eps_{t,x}(y)>0$. Therefore, $\log p_{t,x}^\epsilon(y)$ is well-defined for all $y\in \mathring{\mathcal{S}}$. Notice that $\mathring{\mathcal{S}}\ne \emptyset$.

The next statements provide results on the finiteness of $I$, defined in \eqref{eq:ratefunc}.
\begin{prop}
\label{finite}
The hypotheses are the same as in Theorem \ref{ta.1}. In addition, we suppose that $\sigma,b\in\caC^1$ are Lipschitz continuous and have bounded derivatives. Then, for all $z\in \mathring{\mathcal{S}}$, $I(z)<\infty$.
\begin{proof}
Let $z\in \mathring{\mathcal{S}}$ and $\rho^z:=\dist(z,\partial \mathcal{S})$. Define $z^*_1:=z-\rho^z/2$ and $z^*_2:=z+\rho^z/2$. Since the set $S_1:=\{\Phi_{t,x}^{h};\; h\in \caH_T\}$ is dense in $\mathcal{S}$, there exists $z_1\in S_1\cap B_{\rho^z/4}(z^*_1)$ and $z_2\in S_1\cap B_{\rho^z/4}(z^*_2)$. By definition of $I$ in \eqref{eq:ratefunc}, $I(z_1)$ and $I(z_2)$ are finite.
The function $h\mapsto \Phi^h_{t,x}$ is continuous (see Lemma \ref{lem:frechetdiff}). Hence, by the intermediate value theorem (\cite[Theorem 24.3]{munkres}), we conclude that for all $\bar z\in(z_1,z_2)$ there exists an $h_{\bar z}\in \caH_T$ such that $\Phi^{h_{\bar z}}_{t,x}=\bar z$ and therefore $I(z)<\infty$.
\end{proof}
\end{prop}

If the function $b$ is bounded, one can show that $\{z\in\R;\; I(z) < \infty\}=\R$. Therefore $\supp(P\circ[u^\eps(t,x)]^{-1})=\R$, $p^\eps_{t,x}(y)>0$ for all $y\in\R$ and \eqref{eq:ldpdenslower}
holds for any $y\in\R$. This is a consequence of the following Proposition.

\begin{prop}\label{prop:positivitydensity}
Assume \ref{itm:finiteLambda}, \ref{itm:conditionic}, \ref{itm:sigmabound}, and that $\sigma$ and $b$ are Lipschitz continuous. Suppose also that $b$ is bounded. Then $\{z\in\R;\; I(z) < \infty\}=\R$.
\begin{proof}
Fix $(t,x)\in(0,T]\times\Rd$. Owing to \ref{itm:finiteLambda}, for every $h\in \caH_T$,

\begin{equation}
	\bigg|\int_0^t\int_\Rd \Lambda(t-s,x-z) b(\Phi^{h}_{s,z})\intd z\intd s\bigg| \leq |b|_\infty \int_0^t \Lambda(s)(\Rd) \intd s =: I_2. \label{eq:prf391}
\end{equation}
Note that this bound does not depend on $h\in \caH_T$. Moreover, \ref{itm:sigmabound} and \ref{itm:finiteLambda} imply
\begin{equation}
	I_1 := \sigma_0^2 \|\Lambda(t-\cdot,x-\ast)\|^2_{\caH_T} < \infty.
	\label{eq:prf392}
\end{equation}
Fix $\alpha>0$, $z\in\R$, $h\in \caH_T$, and set
\[ h_{z,\alpha}(\cdot,\ast) := \frac{|z|+\alpha+I_2+|w(t,x)|}{I_1}
\Lambda(t-\cdot,x-\ast)\sigma(\Phi^h_{\cdot,\ast}). \]
Using \eqref{eq:finitebhphi} one can easily check that $h_{z,\alpha}\in \caH_T$.
From \eqref{eq:prf391}, \eqref{eq:prf392}, along with \ref{itm:sigmabound}, we obtain
\begin{align*}
	\Phi^{h_{z,\alpha}}_{t,x}
  & = w(t,x) +  \frac{|z|+\alpha+I_2+|w(t,x)|}{I_1}\\
  &\qquad \times\langle\Lambda(t-\cdot,x-\ast)\sigma(\Phi^{h_{z,\alpha}}_{\cdot,\ast}),\Lambda(t-\cdot,x-\ast)\sigma(\Phi^h_{\cdot,\ast})\rangle_{\caH_T} \notag\\
  & \phantom{=} + \int_0^t\int_\Rd \Lambda(t-s,x-z)b(\Phi^{h_{z,\alpha}}_{s,z})\intd z\intd s \\
  & \geq -|w(t,x)|+\frac{|z|+\alpha+I_2+|w(t,x)|}{I_1} I_1 - I_2 \\
  & > z,
\end{align*}
and similarly,
\begin{align*}
	\Phi^{-h_{z,\alpha}}_{t,x}
  & = w(t,x) -\frac{|z|+\alpha+I_2+|w(t,x)|}{I_1}\\
  &\qquad \times\langle\Lambda(t-\cdot,x-\ast)\sigma(\Phi^{-h_{z,\alpha}}_{\cdot,\ast}),\Lambda(t-\cdot,x-\ast)\sigma(\Phi^h_{\cdot,\ast})\rangle_{\caH_T} \notag\\
  & \phantom{=} + \int_0^t\int_\Rd \Lambda(t-s,x-z) b(\Phi^{-h_{z,\alpha}}_{s,z})\intd z\intd t \\
  & \leq |w(t,x)|-\frac{|z|+\alpha+I_2+|w(t,x)|}{I_1}I_1 + I_2 \\
  & < z.
\end{align*}
Thus, for all $z\in\R$ there exists $h_{z,\alpha}\in \caH_T$ such that $\Phi^{-h_{z,\alpha}}_{t,x} < z < \Phi^{h_{z,\alpha}}_{t,x}$. By the intermediate value theorem \cite[Theorem 24.3]{munkres} together with Lemma \ref{lem:frechetdiff}, there exists some $h_z\in \caH_T$ such that $\Phi^{h_z}_{t,x}=z$. This finishes the proof.
\end{proof}
\end{prop}

Assume as in the previous proposition that \ref{itm:finiteLambda}, \ref{itm:conditionic}, \ref{itm:sigmabound} hold and that $\sigma$ and $b$ are Lipschitz continuous. Suppose that $\sigma$ is bounded. Then, for all $y\in\R$, $p_{t,x}(y)>0$ (see  \cite[Theorem 5.1]{enualart}) and therefore $\mathcal{S}= \R$.


\section{Support theorem}\label{sec:A3}

In this section we prove a characterization of the topological support of the probability law of the random variable $u^\eps(t,x)$, for a fixed $(t,x)\in(0,T]\times\Rd$, defined by \eqref{2.2}. This is the smallest closed subset $\mathcal{X}^\eps\subseteq\R$ satisfying $\P\circ (u^\eps(t,x))^{-1}(\mathcal{X}^\eps)=1$. Under stronger assumptions than in Theorem \ref{basic}, we prove in Theorem \ref{ta.1} that $\mathcal{X}^\eps = \overline{\{\Phi_{t,x}^h; h\in\caH_T\}}$ and therefore also that $\mathcal{X}^\eps$ does not depend on $\eps$.
This will be a consequence of two approximation results, as follows.

Fix $\eps\in(0,1]$ and consider a sequence $(v^{n,\eps})_{n\in\N}$ of $\caH_T$-valued random variables such that
\begin{description}
	\item{(C1)} $\lim_{n\to\infty} \P\big[\big|u^\eps(t,x) - \Phi_{t,x}^{v^{n,\eps}}\big| > \eta\big] = 0$,
	\item{(C2)} for any $h\in\caH_T$, $\lim_{n\to\infty} \P\big[\big|u(t,x; \omega-v^{n,\eps}+h) - \Phi_{t,x}^{h}\big| > \eta\big]= 0$,
\end{description}
for any $\eta>0$.

By Portmanteau's Theorem, (C1) implies that $\mathcal{X}^\eps \subseteq \overline{\{\Phi_{t,x}^h; h\in\caH_T\}}$. From (C2) together with Girsanov's theorem, we deduce the converse inclusion $\mathcal{X}^\eps \supseteq \overline{\{\Phi_{t,x}^h; h\in\caH_T\}}$. Without loss of generality, we may assume that $\eps=1$ and write $u(t,x)$ and $v^n$ instead of $u^\eps(t,x)$ and $v^{n,\eps}$, respectively.

It is easy to see that both convergences (C1) and (C2) (with $\eps=1$) can be formally derived from a single convergence result.
Indeed, let $A,B,G,b:\R\rightarrow\R$, $w:[0,T]\times\Rd\rightarrow\R$ and $h\in\caH_T$. We consider the SPDEs
\begin{align}
	\label{a3}
	X(t,x) = & w(t,x) + \int_0^t\int_\Rd \Lambda(t-s,x-y)(A+B)(X(s,y))M(\intd s,\intd y) \nonumber\\
	& + \langle \Lambda(t-\cdot,x-\ast)G(X(\cdot,\ast)),h\rangle_{\caH_T} \nonumber\\
	& + \int_0^t\int_\Rd \Lambda(t-s,x-y)b(X(s,y))\intd y\intd s,
\intertext{and, for all $n\in\N$}
	\label{a4}
	X_n(t,x) = & w(t,x) + \int_0^t\int_\Rd \Lambda(t-s,x-y)A(X_n(s,y))M(\intd s,\intd y) \nonumber\\
	& + \langle \Lambda(t-\cdot,x-\ast)B(X_n(\cdot,\ast)),v^n\rangle_{\caH_T} \nonumber\\
	& + \langle \Lambda(t-\cdot,x-\ast)G(X_n(\cdot,\ast)),h\rangle_{\caH_T} \nonumber\\
	& + \int_0^t\int_\Rd \Lambda(t-s,x-y)b(X_n(s,y))\intd y\intd s.
\end{align}
Suppose we can prove that the sequence $(X_n(t,x))_{n\in\N}$ converges in probability to $X(t,x)$, for fixed $(t,x)\in(0,T]\times\Rd$. Then with the choice $A=G=0$ and $B=\sigma$, we obtain (C1). By taking $A=G=\sigma$ and $B=-\sigma$, we get (C2).

The sequence $(v^n)_{n\in\N}$ will  consists of smooth approximations of the stochastic process $F$. As has been described in Section \ref{sec:2}, $F$ can be identified with a sequence of independent standard Brownian motions $W=\{W^k(t), t\in[0,T], k\in\N \}$.

Fix $n\in\N$ and consider the partition of $[0,T]$ determined by $\tfrac{iT}{2^n}$, $i=0,1,\ldots,2^n$. Denote by $\Delta_i$ the interval $[\tfrac{iT}{2^n},\tfrac{(i+1)T}{2^n}[$ and by $|\Delta_i|=T2^{-n}$ its length. We write $W_j(\Delta_i)$ for the increment $W_j(\tfrac{(i+1)T}{2^n})-W_j(\tfrac{iT}{2^n})$, $i=0,\ldots,2^n-1$, $j\in\N$. Define differentiable approximations of $(W_j, j\in\N)$ as follows:
\begin{equation*}
	W^n=\left(W_j^n=\int_0^\cdot \dot{W}_j^n(s)ds,\, j\in\N \right),
\end{equation*}
where for $j>n$, $\dot{W}_j^n=0$, and for $1\le j\le n$,
\begin{equation*}
	\dot{W}_j^n(t) =
	\begin{cases}
		\sum_{i=0}^{2^n-2} 2^{n}T^{-1}W_j(\Delta_i)1_{\Delta_{i+1}}(t) & \text{if $t\in[2^{-n}T,T]$,}\\
		&	\\
		0 &	\text{if $t\in[0,2^{-n}T[$.}
	\end{cases}
\end{equation*}
Then, let
\begin{equation*}\label{a2}
	v^n(t,x)=\sum_{j\in \N} \dot{W}_j^n(t) e_j(x).
\end{equation*}

By Theorem \ref{thm:existenceanduniqueness}, Equation \eqref{a3} has a unique random-field solution, and this solution possesses moments of any order, uniformly in $(t,x)$.  That theorem cannot be applied to Equation \eqref{a4}, because the $\caH_T$-valued random variable $v^n$ does not satisfy \eqref{eq:finitehHT}. For this reason (but also for others that will become clear later), we fix a parameter $\theta>0$ and introduce a localization on $\Omega$ defined by
\begin{equation*}\label{loc}
	L_n(t) = \bigg\{ \sup_{1\leq j\leq n}\sup_{0\leq i\leq\lfloor 2^ntT^{-1}-1\rfloor_+} |W_j(\Delta_i)| \leq 2^{n(\theta-1)}\bigg\},
\end{equation*}
where $n\in\mathbb{N}$ and $t\in[0,T]$. Notice that $L_n(t)$ decreases with $t$. Similarly as in \cite[Lemma 2.1]{millet--sanz-sole00}, one can prove that if $\theta>\frac{1}{2}$,
\begin{equation}\label{a5}
	\lim_{n\to\infty} \P\left(L_n(t)^c\right)=0, \  t\in[0,T].
\end{equation}
It is easy to check that
\begin{equation*}\label{a50}
	\Vert v^n(t,\ast)1_{L_n(t)}\Vert_\caH \le Cn2^{n\theta},
\end{equation*}
and for $0\le t<t^{\prime}\le T$,
\begin{equation}\label{a500}
	\Vert v^n1_{L_n(t^{\prime})}1_{[t,t^{\prime}]}\Vert_{\caH_T} \le Cn2^{n\theta} |t-t^{\prime}|^{\frac{1}{2}}.
\end{equation}

On each $L_n(t)$, the assumptions of Theorem \ref{thm:existenceanduniqueness} are satisfied. Thus, by localization we can prove the existence of a unique solution to \eqref{a4}, and that this solution is bounded in probability.


For the formulation of the main result, it is necessary to introduce an additional assumption:
\medskip

\noindent {\bf (A7)}
As in {\bf(A1)}, the mapping $t\mapsto \Lambda(t)$ is a deterministic function with values in the space of non-negative tempered distributions with rapid decrease, and for any $t\in[0,T]$, $\Lambda(t)$ is a non-negative measure. Moreover, there exist $\eta,\delta>0$ such that
\begin{description}
	\item{(i)} $ \int_0^t \int_{\R^d}\vert \caF \Lambda(s)(\xi)\vert^2 \mu(\intd\xi) \intd s \le C t^{\eta}$, for any $t\in(0,T]$,
	\item {(ii)} $\sup_{0\le s\le T} \int_{\R^d}\vert \caF \Lambda(s)(\xi)\vert^2 \mu(\intd\xi) = \sup_{0\le s\le T} J_1(s) <\infty$,
	\item {(iii)} $\int_0^t \Lambda(s)(\R^d) \intd s \le C t^\delta$, for any $t\in(0,T]$.
\end{description}
Clearly, {\bf (A7)} is stronger than {\bf (A1)}.

Let $\Gamma(dx)=|x|^{-\beta} dx$, $\beta\in(0,d\wedge 2)$, and therefore $\mu(\intd\xi)= |\xi|^{-(d-\beta)}\intd\xi$. In Section \ref{sec:applications}, we will see that the fundamental solution to the wave equation with $d=\{1,2,3\}$ satisfies {\bf (A7)}. 
\smallskip

The main result of this section is the following.

\begin{thm}\label{ta.1}
The hypotheses are {\bf (A2)} and {\bf (A7)}. We also suppose that $\sigma$ and $b$ are Lipschitz continuous functions. Let $u^\eps(t,x)$ be the solution to \eqref{2.2} at a given point $(t,x)\in(0,T]\times \R^d$. Then the topological support of the probability law $\P\circ (u^\eps(t,x))^{-1}$ is the closure of the set $\{\Phi_{t,x}^h; h\in\caH_T\}$, where $\Phi_{t,x}^h$ is defined in \eqref{eq:defbhphi}.
\end{thm}

\medskip

By the preceding discussion, the theorem is a corollary of the next Proposition.


\begin{prop}
\label{pa.1}
Assume that $A,B,G,b:\R\rightarrow\R$ are Lipschitz continuous functions and that assumptions {\bf(A2)} and {\bf (A7)} are satisfied. Fix $(t,x)\in(0,T]\times\Rd$ and, in the definition of $L_n(t)$, fix $\theta>\tfrac{1}{2}$ such that $\left(\tfrac{\eta}{2}-\theta+\tfrac{1}{2}\right)\wedge \delta-\theta>0$. Then, for any $p\in[1,\infty)$,
\begin{equation}\label{a6}
	\lim_{n\to\infty}\left\Vert (X_n(t,x)-X(t,x)) 1_{L_n(t)}\right\Vert_p = 0,
\end{equation}
where $\Vert \cdot\Vert_p$ denotes the $L^p(\Omega)$-norm.
\end{prop}

Indeed, owing to \eqref{a5}, the convergence \eqref{a6} yields $\lim_{n\to\infty} X_n(t,x)= X(t,x)$ in probability.

\medskip

The rest of the section is devoted to the proof of Proposition \ref{pa.1}. First we introduce some additional notation.

For any $n\in\N$, $t\in [0,T]$, set
\[ \underline{t}_n = \max\left\{k2^{-n}T, k=0,\ldots,2^n-1: k2^{-n}T < t\right\}, \]
and define $t_n= \max \{\underline{t}_n- 2^{-n}T, 0 \}$. To strengthen the $\caF_t$-measurability properties of $X(t,x)$ and $X_n(t,x)$, we consider stochastic processes defined by a modification of Equations \eqref{a3}, \eqref{a4}, respectively, as follows:

\begin{align}
\label{a7}
	X(t,t_n,x) =
	& w(t,x)+ \int_0^{t_n} \int_{\R^d} \Lambda(t-s,x-y) (A+B)(X(s,y)) M(\intd s,\intd y)\nonumber\\
	& +\langle \Lambda(t-\cdot,x-\ast)G(X(\cdot,\ast))1_{[0,t_n]}(\cdot),h\rangle_{\mathcal{H}_T}\nonumber\\
	&	+\int_0^{t_n} \int_{\R^d} \Lambda(t-s,x-y) b(X(s,y)) \intd y\intd s,
\intertext{and}
\label{a60}
	X_n^-(t,x) =
	& w(t,x)+ \int_0^{t_n} \int_{\R^d} \Lambda(t-s,x-y) A(X_n(s,y)) M(\intd s,\intd y)\nonumber\\
	&	+\langle \Lambda(t-\cdot,x-\ast)B(X_n(\cdot,\ast))1_{[0,t_n]}(\cdot),v^n\rangle_{\mathcal{H}_T}\nonumber\\
	&	+\langle \Lambda(t-\cdot,x-\ast)G(X_n(\cdot,\ast))1_{[0,t_n]}(\cdot),h\rangle_{\mathcal{H}_T}\nonumber\\
	&	+\int_0^{t_n} \int_{\R^d} \Lambda(t-s,x-y) b(X_n(s,y)) \intd y\intd s. \nonumber
\end{align}

\smallskip

In the proof of Proposition \ref{pa.1}, we will use the following facts: for any $p\in[1,\infty)$ and every integer $n\ge 1$,
\begin{description}
	\item{(P1)}
		\begin{equation}\label{s3.91}
			\sup_{(t,x)\in[0,T]\times\R^d} \Vert X(t,x)-X(t,t_n,x)\Vert_p\le C 2^{-n(\frac{\eta}{2}\wedge \delta)}
		\end{equation}
	\item{(P2)}
		\begin{equation}\label{s3.92}
			\sup_{n\in\N} \sup_{(t,x)\in[0,T]\times\R^d} \Vert X(t,t_n,x)\Vert_p\le C
		\end{equation}
	\item {(P3)}
		\begin{equation}\label{s3.94}
			\sup_{(t,x)\in[0,T]\times\R^d} \Vert (X_n(t,x)-X_n^-(t,x)) 1_{L_n(t)}\Vert_p\le Cn2^{-n([(\frac{\eta}{2}-\theta+\frac{1}{2})\wedge \delta]}
		\end{equation}
\end{description}
The estimate \eqref{s3.91} can be easily obtained by adapting the arguments of the proof of \cite[Lemma 4.1]{delgado--sanz-sole012} and applying the assumption {\bf(A7)}. From \eqref{s3.91} and Theorem \ref{thm:existenceanduniqueness} applied to the Equation \eqref{a3}, we obtain \eqref{s3.92}. Finally, \eqref{s3.94} is proved adapting the arguments of the proof of \cite [Lemma 4.2, Lemma 4.3]{delgado--sanz-sole012}, and assuming {\bf(A7)}.

\medskip

\noindent{\it Proof of Proposition \ref{pa.1}}.
\smallskip

\noindent
Using \eqref{a3}, \eqref{a4}, we write the difference $X_n(t,x)-X(t,x)$ grouped into comparable terms and prove their convergence to zero. The main difficulty lies in the convergence of $\langle \Lambda(t-\cdot,x-\ast) B(X_n(\cdot,\ast)),v^n\rangle_{\caH_T}$ to $\int_0^t\int_{\R^d}\Lambda(t-s,x-y)B(X(s,y))M(\intd s,\intd y)$.

We write
\begin{equation*}
	X(t,x)-X_n(t,x)=\sum_{i=1}^{10} U_n^i(t,x),
\end{equation*}
where
\begin{align*}
U_n^1(t,x)	&	= \int_0^t\int_{\R^d} \Lambda(t-s,x-y)\left[A(X(s,y))-A(X_n(s,y))\right] M(\intd s,\intd y),\\
U_n^2(t,x)	&	= \langle \Lambda(t-\cdot,x-\ast)[G(X(\cdot,\ast))-G(X_n(\cdot,\ast))], h \rangle_{\mathcal{H}_T},  \\
U_n^3(t,x)	&	= \int_0^t \int_{\R^d} \Lambda(t-s,x-y)[b(X(s,y))-b(X_n(s,y))] \intd y \intd s, \\
U_n^4(t,x)	& = \int_{t_n}^t\int_\Rd \Lambda(t-s,x-y)B(X(s,y))M(ds,dy), \\
U_n^5(t,x)	& = \langle\Lambda(t-\cdot,x-\ast)[B(X(\cdot,\ast)) - B(X_n(\cdot,\ast))]\Infkt{[t_n,t]}(\cdot),v^n\rangle_{\caH_T}, \\
U_n^6(t,x)	& = -\langle\Lambda(t-\cdot,x-\ast)B(X(\cdot,\ast))\Infkt{[t_n,t]}(\cdot),v^n\rangle_{\caH_T}, \\
U_n^7(t,x)	&	= \int_0^{t_n}\int_{\Rd} \Lambda(t-s,x-y)[B(X(s,y))-B(X^-(s,y))] M(\intd s,\intd y),	\\
U_n^8(t,x)	&	= \int_0^{t_n}\int_{\Rd} \Lambda(t-s,x-y)B(X^-(s,y)) M(\intd s,\intd y) \\
						& \qquad - \langle \Lambda(t-\cdot,x-\ast)B(X^-(\cdot,\ast))\Infkt{[0,t_n]}(\cdot), v^n \rangle_{\mathcal{H}_T},\\
U_n^9(t,x)	&	= \langle \Lambda(t-\cdot,x-\ast)[B(X^-(\cdot,\ast)) - B(X_n^-(\cdot,\ast))]\Infkt{[0,t_n]}(\cdot), v^n \rangle_{\mathcal{H}_T}, \\
U_n^{10}(t,x)	&	= \langle \Lambda(t-\cdot,x-\ast)[B(X_n^-(\cdot,\ast))-B(X_n(\cdot,\ast))]\Infkt{[0,t_n]}(\cdot), v^n \rangle_{\mathcal{H}_T}.
\end{align*}
Here, we have used the abridged notation $X^-(\cdot,\ast)$ for the stochastic process $X(t,t_n,x)$ defined in \eqref{a7}. Notice that although not apparent in this new notation, $X^-(\cdot,\ast)$ does depend on $n$.

Fix $p\in[2,\infty[$. Clearly,
\begin{equation*}
	\E\left(\left| X_n(t,x)-X(t,x)\right|^p1_{L_n(t)}\right) \le C \sum_{i=1}^{10} \E\left(\left|	U_n^i(t,x)\right|^p1_{L_n(t)}\right).
\end{equation*}
We start by analyzing the contribution of $U_n^i(t,x)$, $i=1,2,3$ on the left-hand side of this expression.

Burkholder's and H\"older's inequalities yield
\begin{equation}\label{s3.76}
 	\E\left(\left| U_n^1(t,x)\right|^p1_{L_n(t)}\right)
 	\le C\int_0^t \sup_{y \in \Rd} \E\Big(\Big|X(s,y)-X_n(s,y)\Big|^{p} 1_{L_n(s)}\Big)J_1(t-s) \intd s.
\end{equation}

Schwarz's inequality implies
\begin{align*}
	\E\Big[ & \big| U_n^2(t,x)\big|^p \Infkt{L_n(t)} \Big] \\
	& \le \|h\|_{\mathcal{H}_T}^p \E\left(\left\| \Lambda(t-\cdot,x-\ast)[G(X(\cdot,\ast))-G(X_n(\cdot,\ast))]1_{L_n(t)}\right\|_{\mathcal{H}_T}^2\right)^{p/2}.
\end{align*}
Then, by using H\"older's inequality we obtain
\begin{equation}\label{s3.77}
	\E\left(\left| U_n^2(t,x)\right|^p1_{L_n(t)}\right)
 	\le C \int_0^t \sup_{y \in\Rd} \E\Big(\Big|X(s,y)-X_n(s,y)\Big|^{p} 1_{L_n(s)}\Big) J_1(t-s) \intd s.
\end{equation}

We apply H\"older's inequality to $U_n^3(t,x)$ and obtain
\begin{equation}\label{s3.78}
	\E\left(\left| U_n^3(t,x)\right|^p1_{L_n(t)}\right)
	\le C\int_0^t \sup_{y \in \Rd}\E\Big(\Big|X(s,y)-X_n(s,y)\Big|^{p} 1_{L_n(s)}\Big)J_2(t-s)\intd s.
\end{equation}

Next we consider the terms $U_n^i(t,x)$ for $i=4,5,6$.
Let $i=4$. H\"older's inequalities and assumption {\bf(A7)} yield
\begin{align}\label{U_n^4}
	& \E\left(\left| U_n^4(t,x)\right|^p1_{L_n(t)}\right) \notag\\
	& \quad \le C \bigg(1+\sup_{(t,x)\in[0,T]\times\Rd}\E\big[|X(t,x)|^p\big]\bigg) \bigg(\int_{t_n}^t J_1(t-s)\intd s\bigg)^{p/2} \notag\\
	& \quad \leq C 2^{-pn\eta/2}.
\end{align}

Using H\"older's inequality, \eqref{a500} and assumption {\bf(A7)}, we have
\begin{align}\label{U_n^5}
	& \E\left(\left| U_n^5(t,x)\right|^p1_{L_n(t)}\right) \notag\\
	& \leq Cn^p2^{np(\theta-1/2)}\E\big[\|\Lambda(t-\cdot,x-\ast)[B(X(s,y)) -B(X_n(\cdot,\ast))]\Infkt{[t_n,t]}(\cdot)\Infkt{L_n(t)}\|^p_{\caH_T}\big] \notag\\
	& \leq Cn^p2^{np(\theta-1/2)}\bigg(\int_{t_n}^t J_1(t-s) ds\bigg)^{p/2-1} \notag\\
	& \qquad \times	\int_{t_n}^t\sup_{y\in\Rd}\E\Big(\Big|X(s,y)-X_n(s,y)\Big|^{p} 1_{L_n(s)}\Big)J_1(t-s)ds \notag\\
	& \leq Cn^p2^{-n[p(\frac{\eta}{2}-\theta+\frac{1}{2}) - \eta]} \int_{0}^t\sup_{y\in\Rd}\E\Big(\Big|X(s,y)-X_n(s,y)\Big|^{p} 1_{L_n(s)}\Big)J_1(t-s)ds.\notag\\
	\end{align}
Since $\frac{\eta}{2}-\theta+\frac{1}{2}>0$, for $p>\eta(\frac{\eta}{2}-\theta+\frac{1}{2})^{-1}$, we clearly have $p(\frac{\eta}{2}-\theta+\frac{1}{2}) - \eta>0$.

For $U_n^6(t,x)$, we proceed in a similar manner as for $U_n^5(t,x)$ applying the fact that $X(t,x)$ has uniformly bounded moments of all orders. We obtain
\begin{align}\label{U_n^6}
  & \E\left(\left| U_n^6(t,x)\right|^p1_{L_n(t)}\right) \notag\\
  & \quad \leq Cn^p2^{pn(\theta-1/2)}\bigg(1+\sup_{(t,x)\in[0,T]\times\Rd}\E\big[|X(t,x)|^p\big]\bigg)\bigg(\int_{t_n}^t J_1(t-s)ds\bigg)^{p/2} \notag\\
  & \quad \leq Cn^p2^{-pn(\eta/2-\theta+1/2)}.
\end{align}

Finally, we study $U_n^i(t,x)$, $i=7,8,9,10$.
The arguments based on Burkholder's and  H\"older's inequalities and \eqref{s3.91} give
\begin{align}\label{U_n^7}
	\E\left(\left| U_n^7(t,x)\right|^p1_{L_n(t)}\right)
	&	\le C\int_0^t \sup_{y\in\Rd} \E\Big(\Big|X^-(s,y)-X(s,y)\Big|^{p} 1_{L_n(s)}\Big)J_1(t-s)\intd s \notag\\
	&	\le C 2^{-np[\frac{\eta}{2}\wedge \delta ]}.
\end{align}

In the following, let $\tau_n$ be the operator defined on functions $f:[0,T]\times\Rd\rightarrow\R$ by $\tau_n(f)(t,x):=f((t+2^{-n})\wedge T,x)$. Since $t_n<T-2^{-n}$, the restriction "$\wedge T$" is not active on $t\in[0,t_n]$. Let $\pi_n$ be the projection operator from ${\caH_T}$ onto the Hilbert subspace generated by the set of functions
\begin{equation*}
	\{2^nT^{-1}\Infkt{\Delta_i}(\cdot)\otimes e_k(\ast), i=0,\ldots,2^{n}-1, k=1,\ldots,n\}.
\end{equation*}
Note that $\pi_n\circ\tau_n$ is a uniformly bounded operator in $n\in\N$ and $\pi_n\circ\tau_n$ converges to $I_{\caH_T}$ strongly, where  $I_{\caH_T}$ denotes the identity operator on $\mathcal{H}_T$. Moreover, $\Upsilon_t:=\big(\pi_n \circ \tau_n\big)-I_{\mathcal{H}_T}$ is a contraction operator on $\mathcal{H}_T$.

Since $X_n^-(s,\ast)$, $X^-(s,\ast)$ are $\caF_{s_n}$-measurable random variables, by using the definition of $v^n$ one checks that
\begin{align*}
	U_n^9(t,x) = \int_0^{t_n}\int_{\R^d} (\pi_n\circ\tau_n)\big[ & \Lambda(t-s,x-y) \\
											&	(B(X^-_n(s,y))-B(X^-(s,y)))\big] M(\intd s,\intd y).
\end{align*}
Thus, after having applied Burkholder's inequality, we obtain
\begin{align*}
	& \E\left[\left(\left| U_n^9(t,x)\right|^p\right)1_{L_n(t)}\right]\\
	&	\le C \E\left(\Big\| (\pi_n\circ\tau_n)
\big[\Lambda(t-\cdot,x-\ast)(B(X^-)-B(X_n^-))\big](\cdot,\ast)\Infkt{[0,t_n]}(\cdot)1_{L_n(\cdot)}\Big\|_{\mathcal{H}_T}^p\right)\\
	&	\le C \E\left(\int_0^{t_n} \left\|\Lambda(t-s,x-\ast)(B(X^-(s,\ast))-B(X_n^-(s,\ast)))1_{L_n(s)}
\right\|_{\mathcal{H}}^2 \intd s\right)^{\frac{p}{2}}.
\end{align*}
Then similarly as for  $U_n^2(t,x)$, we have
\begin{align*}
	& \E\left[\left(\left| U_n^9(t,x)\right|^p\right)1_{L_n(t)}\right]	\notag\\
	&	\qquad\le C \int_0^{t_n} \sup_{y\in\Rd}\E\Big(\Big|X^-(s,y)-X_n^-(s,y)\Big|^p 1_{L_n(s)}\Big)J_1(t-s) \intd s.
\end{align*}
This clearly implies
\begin{align*}
	&	\E\left[\left(\left| U_n^9(t,x)\right|^p\right)1_{L_n(t)}\right]\\
	&	\quad\le C \Bigg(\int_0^{t_n} \sup_{y\in\Rd}\E\Big(\Big|X^-(s,y)-X(s,y)\Big|^p 1_{L_n(s)}\Big)J_1(t-s)\intd s\\
	&	\qquad\qquad + \int_0^{t_n} \sup_{y\in\Rd}\E\Big(\Big|X(s,y)-X_n(s,y)\Big|^p 1_{L_n(s)}\Big)J_1(t-s)\intd s\\
	&	\qquad\qquad + \int_0^{t_n} \sup_{y\in\Rd}\E\Big(\Big|X_n(s,y)-X_n^-(s,y)\Big|^p 1_{L_n(s)}\Big)J_1(t-s)\intd s\Bigg).
\end{align*}
Recall that $X^-(s,y)=X(s,s_n,y)$. By applying \eqref{s3.91} and \eqref{s3.94}, we obtain
\begin{align}\label{s3.80}
	& \E\left[\left(\left| U_n^9(t,x)\right|^p\right)1_{L_n(t)}\right]\nonumber\\
	&	\quad\le C \int_0^{t}\sup_{y\in\Rd}\E\Big(\Big|X(s,y)-X_n(s,y)\Big|^p 1_{L_n(s)}\Big)J_1(t-s)\intd s \notag\\
	& \qquad + Cn^p2^{-np[(\frac{\eta}{2}-\theta+\frac{1}{2})\wedge \delta]}.
\end{align}

For the study of  $U_n^{10}(t,x)$, we first apply Schwarz's inequality. Then \eqref{a500} and \eqref{s3.94} yield
\begin{equation}\label{eq:U10}
  \E\left[\left(\left| U_n^{10}(t,x)\right|^p\right)1_{L_n(t)}\right] \leq Cn^{2p}2^{-np[(\frac{\eta}{2}-\theta+\frac{1}{2})\wedge \delta-\theta]}.
\end{equation}

Finally, we consider $U_n^8(t,x)$. We are assuming that $t>0$. Hence, for $n$ big enough, $t_n-2^{-n}>0$ and $t_n+2^{-n}<t$. Define
\begin{align*}
	U_n^{8,1}(t,x)
	&	=\int_0^{t_n}\int_\Rd \pi_n\Big(\Lambda(t-\cdot,x-\ast)B(X^-(\cdot,\ast)) \\
	& \qquad\quad - \tau_n\big[\Lambda(t-\cdot,x-\ast)B(X^-(\cdot,\ast))\Big)(s,y) M(\intd s,\intd y),\\
	U_n^{8,2}(t,x)
	&	=  \int_0^{t_n}\int_\Rd \Big(\Lambda(t-s,x-y)B(X^-(s,y)) \\
	& \qquad\quad - \pi_n\big[\Lambda(t-\cdot,x-\ast)B(X^-(\cdot,\ast))\big]\Big) M(\intd s,\intd y).
 \end{align*}
Clearly,
\begin{equation*}
	U_n^8(t,x)=	U_n^{8,1}(t,x) + U_n^{8,2}(t,x).
\end{equation*}

To facilitate the analysis, we write $U_n^{8,1}(t,x)$ more explicitly, as follows.
\begin{align}\label{s3.82}
	U_n^{8,1}&(t,x)	= \int_0^{t_n}\int_\Rd \Big(\pi_n\big[\Lambda(t-\cdot,x-\ast)B(X^-(\cdot,\ast))\big](s,y) \nonumber\\
	&	\qquad\quad - \pi_n\big[\Lambda(t-\cdot,x-\ast)B\big(X^-\big(\cdot,\ast\big)\big)\big](s+2^{-n},y)\Big)M(\intd s,\intd y).
\end{align}

For the second integral on the right-hand side of \eqref{s3.82} we perform a change of variable $s+2^{-n}\mapsto s$.  Therefore we obtain
\[ \E\left(\left\vert U_n^{8,1}(t,x)\right\vert^p 1_{L_n(t)}\right)\le C\left(V_n^{8,1}(t,x) + V_n^{8,2}(t,x)\right), \]
where
\begin{align*}
  V_n^{8,1}(t,x)
  & := \E\Bigg[\bigg|\int_{t_n}^{t_n+2^{-n}} \int_\Rd \pi_n\big[\Lambda(t-\cdot,x-\ast)B(X^-(\cdot,\ast))\big](s,y)\\
  & \qquad\qquad\qquad \times M(\intd s,\intd y)\bigg|^p1_{L_n(t)}\Bigg], \\
  V_n^{8,2}(t,x) & := \E\Bigg[\bigg|\int_0^{2^{-n}}\int_\Rd \pi_n\big[\Lambda(t-\cdot,x-\ast)B\big(X^-\big(\cdot,\ast\big)\big)\big](s,y) \\
  & \qquad\qquad\qquad \times M(\intd s,\intd y)\bigg|^p1_{L_n(t)}\Bigg].
\end{align*}

By the usual procedure involving Burkholder's and H\"older's inequalities and \eqref{s3.92} we have
\begin{align}\label{Vn81}
  V_n^{8,1}(t,x)
  & \leq C\bigg(1+\sup_{(t,x)\in[0,T]\times\Rd} \E\big[|X^-(t,x)|^p\big]\bigg)\bigg(\int_{t_n}^{t_n+2^{-n}} J_1(t-s)\intd s \bigg)^{\frac{p}{2}} \notag\\
  & \leq C \bigg(1+\sup_{(t,x)\in[0,T]\times\Rd} \E\big[|X^-(t,x)|^p\big]\bigg)\bigg(\int_0^{2^{-n+1}} J_1(s) ds\bigg) ^{\frac{p}{2}}\notag\\
  & \leq C2^{-np\frac{\eta}{2}},
\end{align}
where in the last inequality, after a change of variable, we have applied (i) of assumption {\bf(A7)}.

Notice that for $s\in[0,2^{-n}]$, $X^-(s,y)=X(s,s_n,y)=w(s,y)$. Therefore, condition (ii) in {\bf(A7)} implies,
\begin{align}\label{Vn82}
	V_n^{8,2}(t,x)
	&\le C\left(1+\sup_{(t,x)\in[0,T]\times\R^d}|w(t,x)|^p\right)\left(\int_0^{2^{-n}}J_1(t-s)\intd s\right)^{\frac{p}{2}} \notag\\
	&\le C 2^{-n\frac{p}{2}}.
\end{align}

Thus, by \eqref{Vn81} and \eqref{Vn82} we have proved the convergence
\begin{equation}\label{s3.85}
	\lim_{n\to \infty}  \sup_{(t,x)\in [0,T]\times\Rd} \E\left(\left|U_n^{8,1}(t,x)\right|^p1_{L_n(t)}\right)=0.
\end{equation}


Let us now consider $U_n^{8,2}(t,x)$. After applying Burkholder's inequality we have
\begin{align*}
 	& \E\left(\big| U_n^{8,2}(t,x)\big|^p 1_{L_n(t)}\right) \\
 	& \quad \le C \E\left(\left\|\big(\pi_n-I_{\mathcal{H}_T}\big) \Big[\Lambda(t-\cdot,x-\ast)B(X^-(\cdot,\ast))\Big] 1_{L_n(\cdot)}\Infkt{[0,t_n]}(\cdot)\right\|_{\mathcal{H}_T}^p\right).
\end{align*}
We want to prove that the right-hand side of this inequality tends to zero as $n\to\infty$.

Set
\begin{equation*}
	Z_n(t,x)= \left\|\big(\pi_n-I_{\mathcal{H}_T}\big) \Big[\Lambda(t-\cdot,x-\ast)B(X^-(\cdot,\ast))\Big] 1_{L_n(\cdot)}\Infkt{[0,t_n]}(\cdot) \right\|_{\mathcal{H}_T}.
\end{equation*}
Since $\pi_n$ is a projection on the Hilbert space $\mathcal{H}_T$, the sequence $\{Z_n(t,x), n\ge 1\}$ decreases to zero as $n\to\infty$ and can be bounded from above by $\sup_{n\in\N} \|\Lambda(t-\cdot,x-\ast)B(X^-(\cdot,\ast))\|_{\mathcal{H}_T}$.
Remember that $X^-(s,y)$ stands for $X(s,s_n,y)$, defined in \eqref{a7}, and therefore it depends on $n$.

Assume that
\begin{equation}\label{s3.87}
	\E\left(\sup_n\left\|\Lambda(t-\cdot,x-\ast) B(X^-(\cdot,\ast)) \right\|_{\mathcal{H}_T}^p\right) <\infty,
\end{equation}
for any $p\in[1,\infty)$. Then, by bounded convergence theorem, we can conclude that  $\lim_{n\to\infty}\E[(Z_n(t,x))^p]=0$.

Let us sketch the main arguments for the proof of \eqref{s3.87}. By considering first a convolution in the space variable of $\Lambda(t-\cdot,x-\ast) B(X^-(\cdot,\ast))$ with an approximation of the identity, and then passing to the limit, we prove

\begin{align}
\label{s3.880}
&\E\left(\sup_n \Vert \Lambda(t-\cdot,x-\ast) B(X^-(\cdot,\ast))\Vert_{\mathcal{H}_T}^p\right)\nonumber\\
&\quad\le C \left(1+ \sup_{(t,x)\in[0,T]\times \R^d} E\left(\sup_{n}\left\vert X(t,t_n,x)\right\vert^p\right)\right)\left(\int_0^t J_1(s) \intd s\right)^{\frac{p}{2}}
\end{align}
(see \cite[Proposition 3.3] {nualartquer} for the arguments).

From the definition of $X(t,t_n,x)$, we see that for the second and third terms in \eqref{a7}, the supremum in $n$ can be easily handled, since they are defined pathwise.
For the stochastic integral term, we consider the discrete martingale
\begin{equation*}
\left\{\int_0^{t_n} \int_{\R^d} \Lambda(s_0-s,x-y)(A+B)(X(s,y)) M(\intd s,\intd y), \caF_{t_n}, n\in\mathbb{N}\right\},
\end{equation*}
where $s_0\in]0,T]$ is fixed. By applying first Doob's maximal inequality and then Burkholder's inequality, we obtain
\begin{align*}
&E\left(\sup_{n}\left\vert \int_0^{t_n} \int_{\R^d} \Lambda(s_0-s,x-y)(A+B)(X(s,y)) M(\intd s,\intd y)\right\vert^p\right)\\
&\quad \le C E\left(\left\vert\int_0^t \int_{\R^d} \Lambda(s_0-s,x-y)(A+B)(X(s,y)) M(\intd s,\intd y)\right\vert^p\right)\\
&\quad \le C E\left(\Vert \Lambda(s_0-\cdot, x-\ast)(A+B)(X(\cdot,\ast))\Vert^{\frac{p}{2}}_{\caH_T}\right).
\end{align*}
Finally, we take $s_0 := t$. Using the property  $\sup_{(t,x)\in[0,T]\times \R^d}E\left(\vert X(t,x)\vert^p\right)$, we obtain that the expression \eqref{s3.880} is finite.

Hence, we have proved that
\begin{equation}\label{s3.88}
	\lim_{n\to\infty}E\left(\vert U_n^{8,2}(t,x)\vert^p 1_{L_n(t)}\right) = 0,
\end{equation}
and from \eqref{s3.85} and \eqref{s3.88}, we conclude
\begin{equation}\label{s3.90}
	\lim_{n\to\infty}\E(|U_n^8(t,x)|^p1_{L_n(t)})=0.
\end{equation}


Taking into account \eqref{s3.76} - \eqref{eq:U10} and \eqref{s3.90}, we see that
\begin{align*}
	&	\E\left(\left| X(t,x)-X_n(t,x)\right|^p1_{L_n(t)}\right)\\
	&	\le  C_n +
C\int_0^t \sup_{y\in\Rd} \E\left(\left|X(s,y)-X_n(s,y)\right|^p 1_{L_n(s)}\right)\big(J_1(t-s)+J_2(t-s)\big) ds,
\end{align*}
where $(C_n,n\ge 1)$ is a sequence of real numbers satisfying $\lim_{n\to\infty}C_n=0$. By applying a version of Gronwall's lemma, see \cite[Lemma 15]{dalang}, we finish the proof of the Proposition.
\hfil\cqd

\section{Examples}
\label{sec:applications}

In this section we illustrate Theorem \ref{thm:LDPdensity} with some examples.
\bigskip


\noindent{\bf Stochastic wave equation}
\smallskip

Let $F$ be the Gaussian process introduced in Section \ref{sec:2}. Consider the family of stochastic wave equations indexed by $\eps\in(0,1]$,
\begin{align}
\label{4.1}
  &\bigg(\frac{\partial^2}{\partial t^2} - \Delta_d\bigg)u^\eps(t,x) = \eps\sigma(u^\eps(t,x))\dot{F}(t,x) + b(u^\eps(t,x)),\, (t,x)\in(0,T]\times \R^d,\nonumber\\
  &  u^\eps(0,x) = u_0(x), \quad \frac{\partial u^\eps}{\partial t}(0,x) = u_1(x),
\end{align}
where $\Delta_d$ stands for the $d$-dimensional Laplacian, and $d\in\{1,2,3\}$.

We write \eqref{4.1} in the {\it mild form} \eqref{2.2} with
\begin{equation*}
	w(t,x) = \int_\Rd \Lambda(t,x-y)u_1(y)\intd y + \frac{\partial}{\partial t}\int_\Rd \Lambda(t,x-y)u_0(y)\intd y,
\end{equation*}
where $\Lambda$ is the fundamental solution to the wave equation.

For any $t\in(0,T]$, $\Lambda(t)$ is a measure with support included in $B_t(0)$ (the closed ball of $\R^d$ centered at zero and with radius $t$), and $\Lambda(t)(\R^d)=t$. For example, if $d=3$, it is the uniform surface measure on $\partial B_t(0)$, normalized by the factor $\frac{1}{4\pi t}$.
It is also well-known that the Fourier transform of $\Lambda$ is
$\caF\Lambda(t)(\xi) = \frac{\sin(2\pi t|\xi|)}{2\pi|\xi|}$ for any $d$
(see e.g. \cite[Chapter 5]{Folland76}).

For the sake of illustration, we will assume that the covariance measure of $F$ is given by
$\Gamma(dx)=|x|^{-\beta}\intd x$, with $\beta\in(0,d\wedge 2)$, although a deeper analysis might allow to go beyond this case. Then $\mu (d\xi):=\caF^{-1}(\Gamma)(d\xi)= |\xi|^{-(d-\beta)} d\xi$ and
\begin{equation*}
\int_{\R^d}\vert \caF\Lambda(s)(\xi)\vert^2 \mu(\intd\xi)= \int_{\R^d}\frac{\sin^2(2\pi s|\xi|)}{4\pi^2|\xi|^{d-\beta+2}} \intd\xi.
\end{equation*}
With the change of variable $\xi\mapsto (2\pi s)\xi$, we easily obtain that the last integral is equal to $C s^{2-\beta}$, with $C>0$. Thus,
\begin{equation*}
\int_0^t\int_{\R^d}\vert \caF\Lambda(s)(\xi)\vert^2 \mu(\intd\xi) \intd s = Ct^{3-\beta}, \ t\in[0,T].
\end{equation*}
Consequently, the assumptions \ref{itm:boundsLambda} and {\bf (A7)} hold with $\gamma=\eta = 3-\beta$, $\bar\eta=1$, $\delta=2$ and  $t_0\in(0,T]$.

In this setting, we have the following result on \eqref{4.1}.

\begin{thm}\label{cor:SWE}
	Assume \ref{itm:conditionic}, \ref{itm:sigmabCinf} and \ref{itm:sigmabound}. Then for all $(t,x)\in(0,T]\times\Rd$,
\begin{equation}
	\label{4.2}
	\lim_{\eps\downarrow0} \eps^2\log p^\eps_{t,x}(y) = -I(y),
\end{equation}
for all $y$ in the interior of the support of $u(t,x)$, where $I$ is defined in \eqref{eq:ratefunc}. If in addition, either $b$ or $\sigma$ is bounded, then \eqref{4.2} holds for any $y\in \R$.
\end{thm}

Notice that under the standing hypotheses, Theorems \ref{basic}, \ref{density} and \ref{ta.1} hold. We refer to \cite[Lemma 4.2]{dalangquer} for sufficient conditions on the functions $u_0$, $u_1$ implying \ref{itm:conditionic}.

Next we comment on the validity of assumption \ref{itm:LDPueps}.
The sample paths of the random field solution to Equation \eqref{4.1} belong to the space $\mathcal{C}^\alpha([0,T]\times \Rd)$ of $\alpha$-H\"older continuous functions of degree $\alpha\in (0,\frac{2-\beta}{2})$ (see \cite[Section 2.1]{dalang-sanzsole} for a summary of results and references).
In the present framework, and for spatial dimension $d=3$, a {\it large deviation principle} (LDP) for \eqref{4.1} in the space $\mathcal{C}^\alpha([0,T]\times \Rd)$, with rate function $J\equiv I$, is established in \cite[Theorem 1.1]{ortizsanz} (see also \cite{victorthesis}).  Its proof is carried out following the variational approach of Budhiraja and Dupuis in \cite{budhiraja-dupuis} (see also \cite{dupuis-ellis}). By the classical {\it contraction principle} of LDP (\cite[Theorem 4.2.1]{dembozeitouni93}), this implies \ref{itm:LDPueps}. The proof in \cite[Theorem 1.1]{ortizsanz} also applies to $d\in\{1,2\}$. For $d=2$ and with a different method, F. Chenal \cite{chenal} establishes the same LDP. For $d=1$, the reduced form of the stochastic wave equation driven by space-time white noise is considered in \cite{leandrerusso}, and logarithmic estimates for the density are proved.


\begin{proof}[Proof of Theorem \ref{cor:SWE}]  From the preceding discussion, we infer that the random field solution to the stochastic wave Equation \eqref{4.1} at a fixed point $(t,x)\in (0,T]\times \R^d$ satisfies the assumptions of Theorem \ref{thm:LDPdensity}, and that $J=I$.
\end{proof}
\bigskip


\noindent{\bf Stochastic heat equation}
\smallskip

Consider the family of stochastic heat equations indexed by $\eps\in(0,1]$,
\begin{align}
\label{4.3}
  \bigg(\frac{\partial}{\partial t} - \Delta_d\bigg)u^\eps(t,x) & = \eps\sigma(u^\eps(t,x))\dot{F}(t,x) + b(u^\eps(t,x)), (t,x)\in(0,T]\times \R^d, \nonumber\\
  u^\eps(0,x) & = u_0(x),
\end{align}
where the process $F$ is the same as in the example of the stochastic wave equation, and $d\in\N$.

As in the preceding example,
we interpret \eqref{4.3} in the {\it mild form} \eqref{2.2} with
\begin{equation*}
w(t,x) = \int_\Rd \Lambda(t,x-y)u_0(y)\intd x,
\end{equation*}
and
\begin{equation*}\label{eq:fundamentalsolutionheatequation}
	\Lambda(t,x) = \frac{1}{(4\pi t)^{d/2}}\exp\bigg(-\frac{|x|^2}{4t}\bigg), \ (t,x)\in[0,T]\times\Rd.
\end{equation*}
Hence, $\Lambda(t)(\R^d)=1$ for any $t\in(0,T]$.

Let $\Gamma$ be as in the previous example. Then, using the change of variable $\xi \mapsto \sqrt s \xi$, we have
\begin{equation}
\label{4.4}
\int_0^t \int_{\R^d} \vert\caF\Lambda(s)(\xi)\vert^2 \mu(\intd\xi) ds
= \int_0^t \int_{\R^d}\exp(-4\pi^2s|\xi|^2) |\xi|^{\beta-d} d\xi ds = Ct^{\frac{2-\beta}{2}}.
\end{equation}

Hence, the assumptions {\bf(A1)}, {\bf(A3)} hold with $\gamma=\frac{2-\beta}{2}$ and $\delta=1$.

\begin{thm}\label{thm:SHE}
Assume \ref{itm:conditionic}, \ref{itm:sigmabCinf}, \ref{itm:sigmabound}, \ref{itm:LDPueps}. Then for all $(t,x)\in(t_0,T]\times\Rd$,
\begin{align*}
	&\lim_{\eps\downarrow0} \eps^2\log p^\eps_{t,x}(y) \le -J(y), \\
	&\lim_{\eps\downarrow0} \eps^2\log p^\eps_{t,x}(y) \ge -I(y). \\
\end{align*}
The upper bound holds for all $y\in\R$, while the lower bound (with $I$ defined in \eqref{eq:ratefunc}) holds for any
$y$ in the interior of the support of $u(t,x)$. If in addition $\sigma$ is bounded, then the lower bound holds for any $y\in\R$.
\end{thm}

Suppose that the function $u_0$ is bounded and H{\"o}lder continuous with exponent $\alpha\in(0,1]$. Then Lemma 4.2 in \cite{dalangquer} implies \ref{itm:conditionic}. 

Finally, we give some remarks on the hypothesis \ref{itm:LDPueps}. In the literature, there are several results on large deviations for different types of stochastic heat equations with boundary conditions. For example, \cite{chenalmillet} deals with a heat equation with $d=1$ on a bounded domain with either Neumann or Dirichlet boundary conditions, driven by a space-time white noise. In \cite{marquezsarra2}, the dimension $d$ is arbitrary, the boundary conditions are of Dirichlet type, and the noise is spatially correlated. Additional relevant references are \cite{sowers}, where non-Gaussian noises are considered; \cite{peszat} in the framework of evolution equations; \cite{budhiraja-dupuis} illustrates the variational method on reaction-diffusion equations. In \cite{milletsanz1996}, Varadhan estimates have been obtained for the stochastic heat equation in spatial dimension one with space-time white noise on bounded domains.

We are not aware of any reference on large deviations for Equation \eqref{4.3} in the present setting. Nevertheless, we believe that using a similar approach as in \cite{ortizsanz}, such a result could be proved and that the rate function coincides with $I$. If this intuition is correct, the assumption \ref{itm:LDPueps} of Theorem \ref{thm:SHE} could be removed and we will have an equality like \eqref{4.2}.

\setcounter{section}{0}
\renewcommand{\thesection}{\Alph{section}}
\section{Appendix}

This section is devoted to the proof of some auxiliary results used in the paper. In the first part, we state a theorem on existence and uniqueness of a random field solution to a class of SPDEs which applies to the different types of equations that appear in the paper. In the second part, we prove an estimate on the $\caH_T$-norm of the deterministic Malliavin matrix.

\subsection{A result on existence and uniqueness of solution}
Let $H_1$ and $H_2$ be separable Hilbert spaces. Consider mappings
\begin{equation*}
\tilde A,\tilde B,\tilde G:H_2\times H_1\rightarrow H_1
\end{equation*}
satisfying
\begin{itemize}
\item for all $y, y^\prime\in H_1$,
\begin{align*}
	&\sup_{x\in H_2} \left(\|\tilde A(x,y) - \tilde A(x,y')\|_{H_1} + \|\tilde B(x,y) - \tilde B(x,y')\|_{H_1}\right.\nonumber\\
	&\left.\qquad
	+ \|\tilde G(x,y) - \tilde G(x,y')\|_{H_1}\right)
	\leq C\|y-y'\|_{H_1}.
	\label{eq:Lipschitz}
	\end{align*}
\item There exists $q\in[1,\infty)$, and for all $x\in H_2$,
\begin{equation*}
	\|\tilde A(x,0)\|_{H_1} + \|\tilde B(x,0)\|_{H_1} + \|\tilde G(x,0)\|_{H_1} \leq C\big(1+\|x\|_{H_2}^q\big).
	\label{eq:growthrate}
\end{equation*}
\end{itemize}
Combining these two estimates yields
\begin{equation}\label{eq:lipschitz+growth}
  \|\tilde A(x,y)\|_{H_1} + \|\tilde B(x,y)\|_{H_1} +\|\tilde G(x,y)\|_{H_1} \leq C\big(1 + \|y\|_{H_1} + \|x\|^q_{H_2}\big).
\end{equation}

Let $V=(V(t,x), (t,x)\in[0,T]\times\Rd)$ be a predictable $H_2$-valued stochastic process such that for all $p\in[1,\infty)$
\begin{equation}
	\sup_{(t,x)\in[0,T]\times\Rd} \E\big[\|V(t,x)\|_{H_2}^p\big] < \infty.
	\label{eq:boundedmoment}
\end{equation}
Let $U_0=(U_0(t,x), (t,x)\in[0,T]\times\Rd)$ be a predictable $H_1$-valued stochastic process such that for all $p\in[1,\infty)$
\begin{equation*}\label{eq:U0}
	\sup_{(t,x)\in[0,T]\times\Rd} \E\big[\|U_0(t,x)\|_{H_1}^p\big] < \infty.
\end{equation*}
Let $h$ be an $\caH_T$-valued random variable such that
\begin{equation}
	\sup_{\omega\in\Omega} \|h(\omega)\|_{\caH_T} < \infty.
	\label{eq:finitehHT}
\end{equation}
Consider the equation on $H_1$,
\begin{align}
	U(t,x) = 	& U_0(t,x) + \int_0^t\int_\Rd\Lambda(t-s,x-y)\tilde A(V(s,y),U(s,y))M(\intd s,\intd y) \notag\\
	& + \big\langle \Lambda(t-\cdot,x-\ast)\tilde G(V(\cdot,\ast),U(\cdot,\ast)),h\big\rangle_{\caH_T} \notag\\
        & + \int_0^t\int_\Rd\Lambda(t-s,x-y)\tilde B(V(s,y),U(s,y))\intd y \intd s. \label{eq:generalSPDE}
\end{align}

The following statement is a generalization of \cite[Theorem 4.3]{dalangquer} and \cite[Theorem 6.2]{sanzbook}.
\begin{thm}\label{thm:existenceanduniqueness}
Assume \ref{itm:finiteLambda} and let $\tilde A$, $\tilde B$, $\tilde G$, $V$, $U_0$, $h$, be given as above. There exists a unique predictable stochastic process $(U(t,x), (t,x)\in[0,T]\times\Rd)$ with values in $H_1$ satisfying
\eqref{eq:generalSPDE} for all $(t,x)\in[0,T]\times\R^d$, a.s. This solution satisfies
\begin{equation}
\label{bd} \sup_{(t,x)\in[0,T]\times\Rd} \E\big[\|U(t,x)\|^p_{H_1}\big] < \infty,
\end{equation}
for all $p\in[1,\infty)$ and is continuous in $L^2(\Omega)$.
\begin{proof}
We use the classical approach based on Picard's iterations, as in  \cite[Theorem 6.2]{sanzbook}, and similar ideas as in \cite[Theorem 4.3]{dalangquer}, extended to a Hilbert space setting. In comparison with \cite[Theorem 6.2]{sanzbook}, Equation \eqref{eq:generalSPDE} has the extra term $\big\langle \Lambda(t-\cdot,x-\ast)\tilde G(V(\cdot,\ast),U(\cdot,\ast)),h\big\rangle_{\caH_T}$. We illustrate with an example how to deal with it.

Fix $p\in[1,\infty)$. The Cauchy-Schwarz inequality and \eqref{eq:lipschitz+growth} yield
\begin{align*}
  & \E\Big[\|\big\langle \Lambda(t-\cdot,x-\ast)\tilde G(V(\cdot,\ast),U(\cdot,\ast)),h\big\rangle_{\caH_T}\|^p_{H_1}\Big] \\
  & \leq \sup_{\omega\in\Omega} \|h(\omega)\|^p_{\caH_T} \E\big[\|\Lambda(t-\cdot,x-\ast)\tilde G(V(\cdot,\ast),U(\cdot,\ast))\|^p_{H_1\otimes\caH_T}\big] \\
  & \leq C \int_0^t \bigg(1+\sup_{(r,y)\in[0,s]\times\Rd} \E\big[\|U(r,y)\|^p_{H_1}\big] + \sup_{(r,y)\in[0,s]\times\Rd} \E\big[\|V(r,y)\|^{pq}_{H_1}\big]\bigg) \\
  & \qquad\qquad \times J_1(t-s)\intd s\\
  & \leq C \int_0^t \bigg(1+\sup_{(r,y)\in[0,s]\times\Rd} \E\big[\|U(r,y)\|^p_{H_1}\big]\bigg)J_1(t-s)\intd s,
\end{align*}
where in the second inequality we have applied \eqref{eq:lipschitz+growth}, and in the last one \eqref{eq:boundedmoment}.
We leave it to the reader to complete all the details of the proof.
\end{proof}
\end{thm}

In the preceding sections, the following particular cases of Equation \eqref{eq:generalSPDE} have been considered.

\begin{description}
  \item[(1)] $H_1=H_2=\R$, $h = 0$,  $\tilde A=\eps\sigma$ and $\tilde B=b$ do not depend on the first coordinate, $U_0=w$. Then $U=u^\eps$ (see \eqref{2.2}).
  \item[(2)] $H_1=H_2=\R$, $\tilde A = 0$, $\tilde G=\sigma$ and $\tilde B=b$ do not depend on the first coordinate, $U_0=w$.  Then $U:=\Phi^h$  (see \eqref{eq:defbhphi}).
  \item[(3)] $H_1=\caH_T$, $H_2=\R$, $h = 0$, $\tilde A(x,y)= \eps\sigma'(x)y$, $\tilde B(x,y)=b'(x)y$, $U_0=\eps\break\Lambda(t-\cdot,x-\ast)\sigma(u^\eps(\cdot,\ast))$, $V=u^\eps$. Then $U:=Du^\eps$  (see \eqref{Dueps}).
  \item[(4)] $H_1=\caH_T$, $H_2=\R$, $\tilde A = 0$, $\tilde G(x,y)=\sigma'(x)y$, $\tilde B(x,y)=b'(x)y$,  $U_0=\break\Lambda(t-\cdot,x-\ast)\sigma(\Phi^h_{\cdot,\ast})$, $V=\Phi^h$.  Then $U:=\Xi^h$ (see \eqref{eq:defXi}).
  \item[(5)] $H_1=H_2=\R$,  $\tilde A=\sigma$ does not depend on the second coordinate, $\tilde G(x,y)=\sigma'(x)y$, $\tilde B(x,y)=b'(x)y$, $U_0 = 0$, $V=\Phi^h$. Then $U:=\caN(h)$ (see \eqref{n}).
\end{description}

There is yet another particular case of Equation \eqref{eq:generalSPDE} that we met in the proof of Lemma \ref{lem:lem45}. It is obtained by shifting the sample paths $\omega$ by $\eps^{-1}h$, where $h\in\caH_T$, as it is shown in the next Lemma.

\begin{lem}\label{cor:defuh}
The hypotheses are \ref{itm:finiteLambda}, \ref{itm:conditionic}, that $\sigma,b$ are Lipschitz continuous. Let $h\in H$ and set $u^{\eps,h}(t,x; \omega):=u^\eps(t,x;\omega+\eps^{-1}h)$, ${(t,x)\in[0,T]\times\Rd}$, where $u^\eps$ is the solution to \eqref{2.2}. Then $(u^{\eps,h}(t,x),{(t,x)\in[0,T]\times\Rd})$ satisfies the equation
\begin{align}
  u^{\eps,h}(t,x)
  & =  w(t,x) + \eps\int_0^t\int_\Rd \Lambda(t-s,x-z)\sigma(u^{\eps,h}(s,z))M(\intd s,\intd z) \notag\\
  	& + \int_0^t\int_\Rd \Lambda(t-s,x-z)b(u^{\eps,h}(s,z))\intd z\intd s \notag\\
  	& + \langle\Lambda(t-\cdot,x-\ast)\sigma(u^{\eps,h}(\cdot,\ast)),h\rangle_{\caH_T}. \label{eq:defuh}
\end{align}
\end{lem}

The Lemma relies on the formula
\begin{align*}
\bigg(\int_0^t\int_\Rd \Lambda(s,y) & Y(s,y)M(\intd s,\intd y)\bigg)(\omega+h) = \langle\Lambda(\cdot,\ast)Y(s,y)(\omega+h),h \rangle_{\caH_T} \\
& + \left(\int_0^t\int_\Rd \Lambda(s,y)Y(s,y)(\omega+h)M(\intd s,\intd y)\right)(\omega),
\end{align*}
where $h\in\caH_T$ and $(Y(s,y), (s,y)\in[0,T]\times \R^d)$ is a predictable stochastic processes such that $(\Lambda(s,y) Y(s,y), (s,y)\in[0,T]\times \R^d)$ is integrable with respect to the martingale measure $M$ (see \cite{dalang} and \cite{dalangquer} for details).
This is proved by considering first step processes $g$, given by $g(t,x,\omega) = \Infkt{(a,b]}(t)\Infkt{A}(x)X(\omega)$ for $0\leq a< b\leq T$, $A\in\caB_b(\Rd)$ and $X$ a bounded, $\caF_a$-measurable random variable, and then passing to the limit.


\subsection{Analysis of the deterministic Malliavin matrix}\label{sec:A2}
In this section we derive an assertion similar to \cite[Lemma 8.2]{sanzbook} for the Fr\'echet derivative of the function $\Phi$, defined in \eqref{9}, \eqref{eq:defXi}. 

\begin{prop}
\label{A.3}
	The assumptions are \ref{itm:finiteLambda}, $\sigma,b\in\caC^1$ with bounded Lipschitz continuous derivatives. Then,
	for all $\rho\in[0,t]$ and $t\in[0,T]$,
	\[ \sup_{(r,z)\in[0,t]\times\Rd} \|\bar{D}\Phi^{h}(r,z)\big\|^{2p}_{\caH_{t-\rho,t} }\leq C\big(g_1(\rho)\big)^p. \]
\begin{proof}
Fix $p\in[1,\infty)$, $t\in[0,T]$, $\rho\in[0,t]$ and $(r,y)\in[0,t]\times\R^d$. Recall that  $\bar{D}\Phi^{h}(r,y)$ is an $\caH_T$-valued random variable given by
\begin{align*}
	\bar{D}_{\circ,\bullet}\Phi^{h}(r,y)
	= & \Lambda(r-\circ,y-\bullet)\sigma(\Phi^{h}_{\circ,\bullet}) \\
		& + \big\langle\Lambda(r-\cdot,y-\ast)\sigma'(\Phi^{h}_{\cdot,\ast})\bar{D}_{\circ,\bullet}\Phi^{h}(\cdot,\ast),h\big\rangle_{\caH_T} \notag\\
  & + \int_0^r\int_\Rd \Lambda(r-s,y-z) b'(\Phi^{h}_{s,z})\bar{D}_{\circ,\bullet}\Phi^{h}(s,z) \intd z\intd s
\end{align*}
(see \eqref{9} and \eqref{eq:defXi}).

We analyze each one of the three terms on the right-hand side of this equation separately.

For the first term, we have
\begin{align*}
  \big\|\Lambda & (r-\circ,y-\bullet)\sigma(\Phi^h_{\circ,\bullet})\big\|^{2p}_{\caH_{t-\rho,t}} \\
  & \leq \bigg(\int_{t-\rho}^r J_1(r-s)\intd s\bigg)^{p-1}\int_{t-\rho}^r \sup_{(v,z)\in[0,s]\times\Rd} \E\big[|\sigma(\Phi^h_{v,z})|^{2p}\big] J_1(r-s)\intd s \\
  & \leq C(g_1(\rho))^p,
\end{align*}
where in the last inequality we have used that $\sigma$ is Lipschitz continuous and also that for each $h\in\caH_T$, the function
$(t,x)\in[0,T]\times \R^d\mapsto \Phi^h(t,x)$ is uniformly bounded. Indeed, this is a consequence of \eqref{bd}, since $\Phi^h(t,x)$
is deterministic.

For the second term, we apply first Schwarz' inequality and then H\"older's inequality. Using that $\sigma'$ is bounded, we obtain
\begin{align*}
  & \big\|\langle  \Lambda(r-\cdot,y-\ast)\sigma'(\Phi^h(\cdot,\ast))\bar{D}_{\circ,\bullet}\Phi^h(\cdot,\ast),h\rangle_{\caH_T}\|_{\caH_{t-\rho,t}}^{2p} \\
  & \quad \leq C\|h\|^{2p}_{\caH_T} \int_0^r \sup_{(v,z)\in[0,s]\times\Rd} \|\bar{D}_{\circ,\bullet}\Phi^h(v,z)\|^{2p}_{\caH_{t-\rho,t}} J_1(r-s)\intd s.
\end{align*}

Finally, for the last term we apply H\"older's inequality with respect to the finite measure $\Lambda(r-s,x-z)\intd z\intd s$ along with the boundedness of $b'$. We obtain,
\begin{align*}
	\bigg\|\int_0^r\int_\Rd & \Lambda(r-s,y-z)b'(\Phi^h_{s,y})\bar{D}_{\circ,\bullet}\Phi^h(s,z)\intd z\intd s\bigg\|^{2p}_{\caH_{t-\rho,t}} \\
	& \quad \leq C\int_0^r \sup_{(v,z)\in[0,s]\times\Rd} \|\bar{D}_{\circ,\bullet}\Phi^h(v,z)\|_{\caH_{t-\rho,t}}^{2p} J_2(r-s)\intd s.
\end{align*}

By applying Gronwall's lemma to the real function
$$s\mapsto \sup_{(r,z)\in[0,s]\times\Rd} \|\bar{D}\Phi_{\circ,\bullet}^h(r,z)\|^{2p}_{\caH_{t-\rho,t}},$$
 we have
\begin{equation*}
 \sup_{(r,z)\in[0,t]\times\Rd} \|\bar{D}_{\circ,\bullet}\Phi^h(r,z)\|^{2p}_{\caH_{t-\rho,t}} \leq C\bigg(\int_0^\rho J_1(s)\intd s\bigg)^p
 =C\left(g_1(\rho)\right)^p,
 \end{equation*}
for all $t\in[0,T]$. This yields the assertion.
\end{proof}
\end{prop}

\medskip


\end{document}